\documentclass[12pt]{article}
\usepackage{JASA_manu}

\usepackage{setspace}

\usepackage[latin1]{inputenc}
\usepackage{amsmath}
\usepackage{amsfonts}
\usepackage{amsthm}
\usepackage{amssymb}
\usepackage[usenames,dvipsnames]{color}

\usepackage{rotating}
\usepackage{enumerate,rotate}
\usepackage{latexsym}
\usepackage{graphicx, graphics}
\usepackage[english]{babel}
\usepackage{natbib}
 \bibpunct{(}{)}{;}{a}{,}{,}
\usepackage{bm}
\usepackage{caption,subcaption}

\def\CF { \mathcal{F}}

\def\b {\beta}
\def\a {\alpha}

\def \th{\theta}
\def\s {\sigma}

\def\a {\alpha}

\def\J {\mathbb{I}}

\def\N {\mathbb{N}}

\makeatletter
\def\blfootnote{\xdef\@thefnmark{}\@footnotetext}
\makeatother

\def \BGOS {\text{Beta-GOS}}

\newtheorem{definition}{Definition} 

\newtheorem{lemma}[definition]{Lemma}
\newtheorem{proposition}[definition]{Proposition}

\newcommand{\Be}{\text{Beta}}

\newcommand{\IG}{\text{Inv-Gamma}}

\usepackage[usenames]{color}
\newcommand{\ech}{\color{black}\rm  }    

\title{Generalized species sampling priors with latent Beta reinforcements}

\author{
\centering \ 
Edoardo M. Airoldi\and Thiago Costa  \and Federico Bassetti \and Fabrizio Leisen \and Michele Guindani\footnote{ Edoardo M. Airoldi is an Associate Professor in the Department of Statistics at Harvard University and an Alfred P. Sloan Research Fellow (airoldi@fas.harvard.edu), Thiago Costa is a PhD candidate in the School of Engineering and Applied Sciences at Harvard University (tcosta@fas.harvard.edu),   Federico Bassetti is an Assistant Professor in the Department of Mathematics at Universita di Pavia, Italy (federico.bassetti@unipv.it),  Fabrizio Leisen is a Senior Lecturer in the Department of  Statistics at University of Kent, UK (fabrizio.leisen@gmail.com), Michele Guindani is an Assistant Professor in the Department of Biostatistics at MD Anderson Cancer Center (mguindani@mdanderson.org). All authors have equally contributed to the manuscript.}
}
\date{}

\begin{document}
\maketitle

\abstract{Many popular Bayesian nonparametric priors can be characterized in terms of  exchangeable species sampling sequences. However, in some applications, exchangeability may not be appropriate. We introduce a {novel and probabilistically coherent family of non-exchangeable species sampling 
sequences characterized by a tractable predictive probability function with weights driven by a sequence of independent Beta random variables. We compare their theoretical clustering properties with those of the Dirichlet Process and the two parameters Poisson-Dirichlet process. The proposed construction provides a complete characterization of  the joint process, differently from existing work.  We then propose the use of such process as  prior distribution in a hierarchical Bayes modeling framework, and we describe a Markov Chain Monte Carlo sampler for posterior inference. We evaluate the performance of the prior and the robustness of the resulting inference in a simulation study, providing a comparison with popular Dirichlet Processes mixtures and Hidden Markov Models. Finally, we develop an application to the detection of chromosomal aberrations in breast cancer by leveraging array CGH data. \newline}

\noindent
 AMS CLASSIFICATION : Primary 62C10; secondary 62G57





\section{Introduction}
Bayesian nonparametric priors
have become increasingly popular  in applied statistical modeling in the last few years. Examples of their wide area of applications range from variable selection in genetics \citep{Kim06} to linguistics \citep{Teh06b,Wallach08},  psychology \citep{Navarro06}, human learning \citep{Griffiths07},  image segmentation \citep{Sudderth09} and applications to the neurosciences \citep{Jbabdi09}. See also \cite{Hjort10}.
The increased interest in non-parametric Bayesian approaches is motivated by a number of attractive inferential properties. 
For example, Bayesian nonparametric priors are often used as flexible models to describe the heterogeneity of the population of interest, as they implicitly induce a clustering of the observations into homogeneous groups.  Such a clustering can be seen as a realization of a random partition
scheme and can often be characterized in terms
 of a species sampling
 (SS) allocation rule. More formally, a SS sequence  is a  sequence of random variables  $X_1, X_2, \ldots,$
characterized by the predictive probability functions,
\begin{equation}\label{eq:SS_Pit}
P\{ X_{n+1}  \in \cdot \,\,| X_1, \ldots, X_n \} =
\sum_{j=1}^{n} q_{n,j} \delta_{X_{j}}(\cdot)+q_{n,n+1} G_0(\cdot),
\end{equation}
where $\delta_{x}(\cdot)$ denotes a point mass at $x$,   $q_{n,j}$  ($j=1,\dots,n+1$) are non--negative functions of
$(X_1,\dots,X_n)$, or weights, such that $\sum_{j=1}^{n+1} q_{n,j}=1$, and $G_0$
is a non-atomic probability measure \citep{Pit96b}. Collecting the unique values of $X_j$,  \eqref{eq:SS_Pit} can be rewritten as
\begin{equation}\label{eq:SS_Pit2}
P\{ X_{n+1}  \in \cdot \,\,| X_1, \ldots, X_n \} = \sum_{j=1}^{K_n} q_{j}^*\delta_{X_{j}^*}(\cdot)+
q_{K_n+1}^* G_0(\cdot),
\end{equation}
where $K_n$ is the (random) number of distinct values, say $(X_1^*, \ldots X_{K_n}^*)$,
in the vector $X(n)=(X_1, \ldots, X_n)$ and $ q_{j}^*$ are suitable non--negative weights.
In particular, an exchangeable SS sequence is characterized by weights  $ q_{j}^*$  that depend
only on $\mathbf{n}_n=(n_{1n}, \ldots, n_{K_nn})$, where  $n_{jn}$
is the frequency of  $X_j^*$ in $X(n)$ \citep{Fortini2000,han-pitman, Lee08}.
The most well known example of predictive rules of type \eqref{eq:SS_Pit}
is the Blackwell MacQueen sampling rule, which implicitly defines a Dirichlet Process  (DP, \citealp{black-mac,IshZar03}).
The predictive rule characterizing a DP with mass parameter
$\theta$ and base measure $G_0(\cdot)$, $DP(\theta,G_0)$,
sets $q_{n,j}=\frac{1}{n+\theta}$
and $q_{n,n+1}=\frac{\theta}{n+\theta}$ in  \eqref{eq:SS_Pit}.

Whenever the weights $q^*_{j}(\mathbf{n}_n)$ and $q^*_{K_n+1} (\mathbf{n}_n)$ do not depend only on $\mathbf{n}_n$, the sequence $(X_1, X_2,\ldots )$ is not exchangeable. Models  with non-exchangeable random partitions have recently appeared in the literature, e.g.
 to allow for partitions that depend on covariates.
\cite{ParkDunson07} derive a generalized product partition model (GPPM) in which the partition process is predictor--dependent. Their GPPM generalizes the DP clustering mechanism to relax the exchangeability assumption through the incorporation of predictors, implicitly defining a generalized P\'olya urn scheme. \cite{Mueller10} define a product partition model that includes a regression on covariates, allowing units with similar covariates to have  greater probability of being clustered together. Arguably, the previous models provide an implicit modification of the predictive rule \eqref{eq:SS_Pit} where the weights can be  seen as function of some external predictor. Alternatively, other authors model the weights $q_j({\mathbf{n}_n})$ explicitly, for instance, by specifying the weights  as a function of the distance between data points \citep{Dhal08,Blei09}. However, the general properties of the random partitions generated by such processes have not been specifically addressed.

In this paper, we introduce a novel and probabilistically coherent family of non-exchangeable species sampling sequences, where the weights are specified sequentially and do not depend on the cluster sizes, but instead they depend on the realizations of a set of latent variables.
Working within this family, we propose a simple characterization of the weights in the predictive probability function  as a product of independent Beta random variables. This strategy leads to a well-defined random allocation scheme for the observables. The resulting process, which we call Beta-GOS process, is a special case of a Generalized Ottawa Sequence (GOS), recently introduced by \cite{bas-cri-lei}. 

In Section \ref{sec:gos-prior}, we discuss the properties of the Beta-GOS process, with particular regard to the  clustering induced on the observables. In Section \ref{s:priorproperties}, we study the asymptotic distribution of  the (random) number of distinct  values in the sequence, $K_n$, for some natural specifications of the weights, and we compare those results with the well-known asymptotic results characterizing the  DP and the two-parameters Poisson Dirichlet process.
In many applications,  nonparametric processes are often used within hierarchical models to specify the prior distribution of some parameters of the distribution of the observables. For example, this is a popular use for mixtures of Dirichlet Processes. Similarly, the Beta-GOS process can also be used to define a prior  in a hierarchical model. In Section \ref{sec:gos-model}, we outline a basic hierarchical model based on the Beta-GOS process and we outline the basic steps of a MCMC sampler for posterior inference. 
In Section \ref{sec:sim}, we design a set of simulations to we compare the behavior of the Beta-GOS model with that of DP mixtures and hidden Markov Models (HMM) in terms of cluster estimation. Our results suggest that the Beta-GOS  can be seen as a robust alternative to the Dirichlet process when exchangeability would be hardly justified in practice, but still there is a need to describe the heterogeneity of our observations by virtue of an unsupervised clustering of the data. The Beta-GOS  also provides an alternative to customary HMM, especially when the number of states is unknown or the 
Markovian structure is expected to vary with time. 

In Section \ref{sec:app-aberrations}, we analyze two published data sets of genomic and transcriptional aberrations \citep{chin:2006, Curtis12}. Bayesian models for Array CGH data have been recently investigated by \cite{Guha08}, \cite{DeSantis09}, \cite{Baladandayuthapani10}, \cite{Du10}, \cite{Cardin11}, and \cite{Yau11}, among others. Guha et al. propose a four state homogenous Bayesian HMM  to detect copy number amplifications and deletions and partition tumor DNA into regions (clones) of relatively stable copy number. DeSantis et al. extend this approach and develop a supervised Bayesian latent class approach for classification of the clones that relies on a heterogenous hidden Markov model to account for local dependence in the intensity ratios. In a heterogeneous hidden Markov model, the transition probabilities between states depend on each single clone or the the distance between adjacent clones \citep{Marioni06}. Du et al. propose a sticky Hierarchical DP-HMM \citep{Fox11, Teh06a} to infer the number of states in an HMM, while also imposing state persistence. \cite{Yau11} also propose a nonparametric Bayesian HMM, but use instead a DP mixture to model the likelihood in each state. With respect to those proposals, we also assume that the number of states is unknown, as it is typical in a Bayesian nonparametric setting, but we don't need  a parameter to explicitly account for state persistence. This is because the Beta-GOS model is ``non-homogenous'' by design, as the weights in the species sampling mechanisms adapt to take into account the local dependence in the clones'  intensities. We show that the Beta-GOS is able to identify clones that have been linked to breast cancer pathophysiologies in the medical literature. 

We conclude with some final remarks in Section \ref{sec:last}. Technical details and proofs of theorems and lemmas are provided in the Appendix.


\section{The Beta-GOS process.}
\label{sec:gos-prior}

As anticipated, the Beta-GOS process is defined by a modification of the predictive rule that characterizes the species sampling mechanism  \eqref{eq:SS_Pit}, where the weights are a product of independent Beta random variables.
More in general, we start considering a sequence of random variables $(X_n)_{n \geq 1}$ characterized by
the predictive distributions
\begin{equation} \label{betagos1}
   P\{ X_{n+1} \in \cdot \,\, |X(n), W(n)\}=\sum_{j=1}^{n} p_{n,j} \delta_{X_j}(\cdot)+
   r_n G_0(\cdot)
  \end{equation}
where $W(n)=(W_1,\dots,W_n)$ is a vector of independent random variables
$W_k$  taking values in $[0,1]$, and
the weights are defined by
\begin{equation}\label{formapesi}
p_{n,j}:=(1-W_j) \prod_{i=j+1}^{n}W_i, \qquad r_n:=\prod_{i=1}^n W_i.
\end{equation}
The prediction rule
\eqref{betagos1} defines a special case of a Generalized Ottawa Sequence, introduced in
 \cite{bas-cri-lei}, a type of Generalized P\'olya Urn sequences  where the reinforcement is randomly determined by the realizations of a latent process \citep[see also][for an alternative proposal]{Guha10}.
Except from a few special cases, the $X_i$'s in a GOS are not exchangeable. However, it can be shown that these sequences maintain some properties typical of exchangeable sequences. Most notably, any GOS is
 {\it conditionally identically distributed} (CID), i.e.
for all $n> 0$, the $X_{n+j}$'s, $j \geq 1$, are identically
distributed, conditionally on $(X_1, \ldots, X_n, W_1, \ldots,
W_n)$. Hence, the $X_i$'s are also marginally identically
distributed. Note that a CID sequence  is  not necessarily
stationary. If a CID sequence is also stationary then it is
exchangeable. Finally, although no representation theorem is known for CID sequences, it can be shown that given any bounded and measurable function $f$,
the predictive mean $E[f(X_{n+1})|X_1, . . . , X_n]$ and the empirical mean
$\frac{1}{n}\sum_{i=1}^n X_i$  converge to the same limit
as $n$ goes to infinity. For details, see  \cite{be-pra-ri}, where CID sequences have been first introduced. 
The predictive rule \eqref{betagos1} reduces to known cases with a suitable choice of the latent $W_n$'s; for instance  if $W_n:=(\theta+n-1)/(\theta+n)$, then
\eqref{betagos1}  coincides with the Blackwell-MacQueen sampling rule characterizing a
$DP(\theta,G_0)$.

In this paper,  we propose $(W_n)_{n \geq 1}$  be a sequence of independent Beta$(\alpha_n, \beta_n)$ random variables and we call the resulting $(X_1, X_2, \ldots)$ a Beta-GOS sequence.  The choice of Beta latent variables allows for a flexible specification of the species sampling weights, while retaining a simple and interpretable model together with computational simplicity (see later Sections). 
The allocation rule can also be described in terms of a preferential attachment scheme, where each observation is attached to any of the preceding by means of a ``geometric-type'' assignment.
In this scheme, every individual $X_i$ is characterized by a random weight (or ``mark''), $1-W_i$. We can interpret each individual mark as an individual specific attractivity index, as it determines the probability that the next observation will be clustered with $X_i$. More precisely, the first individual is assigned a random value (or ``tag'') $X_1$, according to $G_0$. Now, suppose we have $X_1, \dots, X_n$ together with their marks up to time $n$, $(1-W_1,\dots,1-W_n)$. Then, the $(n+1)$-th individual will be assigned the same tag as $X_n$ with probability $1-W_n$; the probability of pairing $X_{n+1}$ to $X_{n-1}$ will be $W_n (1-W_{n-1})$,  and so forth. In general,  $p_{n,j}$ will be the  product of the repulsions $W_i$ for  the latest $n-j$ subjects and the $j$th attractivity $1-W_j$. Summarizing, $X_{n+1}$ will result in a new tag (i.e., $X_{n+1} \sim G_0$) with probability $r_n$, or will be clustered together with a previously observed tag, say $X_k^*$,  with probability $\sum_{j: X_j=X_k^*}\,  p_{n,j} $.  
In the next Section, we discuss  the clustering behavior induced by different specifications of the Beta weights in more detail.

\section{Clustering behavior of the Beta-GOS.}
\label{s:priorproperties}

The predictive rule \eqref{betagos1} implicitly defines a random partition of the set $\{1,\ldots, n\}$
into $K_n$ blocks. In probability theory,  $K_n$ is also referred to as the {\it  length} of the partition.  Knowledge of the behavior
of  $K_n$ is useful to understand the clustering structure implied by \eqref{betagos1}.
For instance, for a $DP(\theta, G_0)$, it is well-known that $K_n/\log(n)$ converges almost surely to a constant, indeed the mass parameter $\theta$.  This asymptotic behavior is sometimes described as a ``self-averaging" property of the partition \citep{Aoki08}. From a practical point of view, since $K_n/\log(n)$ converges to a constant, then in the limit $K_n$ is essentially $\th \log(n)$; thus, for modeling purposes it suffices to consider only the first moment of $K_n$.
In the case of the two parameter Poisson Dirichlet process
 the length of the partition $K_n$ (suitably rescaled)  converges instead to a random variable.
More precisely, for a $PD(\alpha, \theta)$, with
$0<\alpha<1$, $\theta>-\alpha$,
 then $   K_n/n^{\alpha} $ converges a.s. to  a   strictly positive random variable
(see Theorem 3.8 in \citealp{Pit06}).
Therefore, the PD sequence is non self-averaging.
When the limit of $K_n$ is essentially a random variable, extra care is needed in the prior assessment of the parameters of the nonparametric prior, since the clustering behavior is ultimately governed by the whole distribution of the limit random variable. 
For the Beta-GOS process, we  focus  on the following two cases: 
\begin{itemize}
\item[(i)] $\a_n=a>0$ and $\b_n=b>0$ for all $n\geq 1$;
\item[(ii)] $\a_n=\theta-1+n$ ($\theta>0$) and $\b_n\geq 1$ for all $n\geq 1$ .
\end{itemize}
Then, we can prove the following
\begin{proposition}\label{prop-Kn}
Let $K_n$ be the length of the partition induced by a Beta-GOS, with $G_0$ non-atomic and
$W_n \sim \Be(\alpha_n,\beta_n)$ ($n\geq 1$).
\begin{itemize}
\item[(a)] If $\a_n=n+\th-1,\b_n=1$,  for given $\th>0$,
\( K_n/\log(n)\)
 converges in distribution to a $Gamma(\th,1)$ random  variable.
 \item[(b)] If $\a_n=n+\th-1,\b_n=\beta$, $(\th>0,\beta>1)$
 or if $\a_n=a,\b_n=b$, $(a>0, b>0)$, then $K_n$ converges almost surely
 to a finite random variable $K_\infty$.
In particular, if $\a_n=a,\b_n=b$, then
 \[
E[e^{-t K_\infty}] =e^{-t} \sum_{m \geq 0} (e^{-t}-1)^m \prod_{j=1}^m \frac{(a)^{(j)}}{(a+b)^{(j)}-(a)^{(j)}}
\]
where 
$(t)^{(m)}=t(t+1)\dots(t+m-1)$ and
\[
E[(K_\infty-1)\dots(K_\infty-m)]= m!  \prod_{j=1}^m \frac{(a)^{(j)}}{(a+b)^{(j)}-(a)^{(j)}}.
\]
 \end{itemize}
\end{proposition}

\noindent The proof is detailed in the Appendix,
 where we also provide a general formula for the probability distribution,
 the $k$-th moment and the generating function of $K_n$.
The result in Proposition \ref{prop-Kn}(a) represents a case of a quite natural (non exchangeable) partition model for which the     length $K_n$ scale as $\log(n)$ but is  not self-averaging. When $\a_n=a,\b_n=b$, according to Proposition \ref{prop-Kn}(b),
the convergence of  $K_n$ to a finite
random variable naturally implies the creation of a few big clusters, as $n$ increases.
Instead,  for $\a_n=n+\theta-1,\beta_n=1$, the mean length of the partition depends on the
value of $\theta$, since a bigger number of clusters is associated on average with greater values of $\theta$. However,  as $\theta$ increases so does the asymptotic variability of $K_n$; therefore, in this case, a Beta-GOS  process can be used to represent uncertainty on $K_n$ (by the lack of the self-averaging property of the process. By means of simulations, we have also confirmed that, for small values of $\theta$, the  partition of $n$ elements is skewed, i.e. it is characterized by a small number of big clusters as well as a few small clusters. As $\theta$ increases, the sizes of the clusters decrease accordingly, the observations being grouped into clusters of relatively fewer elements. This is similar to what happens for the DP, and indeed in this case the parameter $\theta$ could be interpreted as a mass parameter for the Beta-GOS.

The parameters of the  $W_i$'s can be chosen to model the autocorrelation expected {\it a priori} in the dynamics of the sequence. The probability of a tie may decrease with $n$ and atoms that have been observed at farthest times may have a greater probability to be selected if they have also been observed more recently. Such considerations may be helpful to guide  prior assessment  of the Beta hyper-parameters.  For given  $n \geq 1$, taking expectations with respect to the weights $W_i$'s we obtain
\begin{equation}\label{Eweights-1}
 E[r_n]=\prod_{j=1}^n \frac{\a_j}{\a_j+\b_j}, \qquad
 E[p_{n,k}]= \frac{\b_k}{\a_k+\b_k} \prod_{j=k+1}^n \frac{\a_j}{\a_j+\b_j}  \qquad k=1,\dots,n.
\end{equation}

Under (a),  it follows that $E[r_n]=({a}/{(a+b)})^n$ and
$E[p_{n,k}]= ({a}/{(a+b)})^{n-k} ({b}/{(a+b)})$; hence, the probabilities of ties depend only on the lag $n-k$ and decrease exponentially as a function of $n-k$. Under (b),
\[
\begin{split}
& E[r_n] = \prod_{j=1}^n \frac{j+\theta-1}{j+\theta-1+\beta}=\frac{\Gamma(\theta+n)\Gamma(\theta+\beta)}{\Gamma(\theta+\beta+n)
\Gamma(\theta)} \\
&  E[p_{n,k}]= \frac{\beta}{k+\theta-1+\beta}  \prod_{j=k+1}^n \frac{j+\theta-1}{j+\theta-1+\beta}
=\beta \frac{\Gamma(\theta+n)\Gamma(\theta-1+\beta+k)}{\Gamma(\theta+\beta+n)
\Gamma(\theta+k)}
\quad k=1,\dots,n. \\
\end{split}
\]
Thus, for $n,k \to +\infty$, $E[r_n] \sim \frac{1}{n^{\beta}}$ and $E[p_{n,k}] \sim \frac{k^{\beta-1}}{n^{\beta}}$. For example, if $\theta=1$ and $\beta=2$, then $\a_j=j$ and $\b_j=2$ and $E[r_n]= \frac{2}{(1+n)(2+n)}$, $E[p_{n,k}]=   \frac{2(k+1)}{(n+1)(n+2)}$, $k=1,\dots,n$, so that the weights decrease linearly as a function of the lag $n-k$.
If  $\alpha_j=\theta-1+j$ ($\theta>0$) and $\beta_j=1$ then $E[r_n]=  \frac{\theta}{\theta+n}$ and $E[p_{n,k}]=   \frac{1}{\theta+n}, k=1,\dots,n$, i.e. any  observation has the same weight. This latter specification leads to an expression similar to that in the Blackwell-McQueen P\'olya Urn characterization of the Dirichlet process; however, this identity is true only in expectation, and the clustering behavior of the DP and 
Beta-GOS  process with $\a_j=j+\theta-1$ and $\b_j=1$ may be quite different, as it is evident from Proposition \ref{prop-Kn}.

In practice, the determination of the parameters of the Beta
distributions is not trivial, and may be problem
dependent, especially given the sensitivity of the clustering behavior to the values of $\alpha_j$ and $\beta_j$. As a general rule, following what it is usually done with Dirichlet processes priors, one should consider eliciting the parameters on the basis of the expected number of clusters $E(K_n)=1+\sum_{j=1}^{n-1}E[r_j]$. For example,    one should set  $\alpha_j=a$ and $\beta_j=b$ to represent a short memory
process,  and the values of $a,b$ can be chosen based on the asymptotic relationship $E(K_n)\approx\frac{a+b}{b}$.  We further suggest to choose $b=1$, or anyway $b<a$, to encourage  a priori low autocorrelation of the sequence, since then $E(p_{n,n})<0.5$. As a matter of fact, we implemented those suggestions in the application to the array CGH data presented in Section \ref{sec:app-aberrations}, where biological considerations lead to further expect the true number of states to be around 4.  
On the other hand,  one should set $\a_j=j+\th-1$, $\b_j=1$ to represent a long memory process, and then choose $\th$ based on 
$E(K_n)=\sum_{j=0}^{n-1}\frac{\theta}{\theta+j} \sim \theta \log(n)$,
for large $n$. The latter, single-parameter, formulation should be  the default choice in those
applications where prior information on the expected number of clusters is unavailable otherwise. 
Alternative strategies are possible. For example, one could consider second moments, or otherwise require further constraints on the expected autocorrelation of the sequence. However, we leave the exploration of those possibilities  to future work. See also the discussion at the end of Section \ref{sec:MCMCsampling}.

Finally, we note the functional form of \eqref{formapesi} may initially suggest a relationship with the stick-breaking characterization of the Dirichlet process.  However, the stick-breaking construction characterizes the representation of the DP as a random measure, not the corresponding predictive probability function. Furthermore, the sequence generated by a DP is exchangeable, whereas a Beta-GOS in general is not and includes the DP as a special case. As a matter of fact, if one would like  to stress the ``stick-breaking''  analogy anyway, one should more properly interpret  \eqref{betagos1} in terms of an {\it inverse} stick-breaking, since each $p_{n,j}$,  which defines the probability of a tie, say $X_{n+1}=X_j$, does not depend on the  $W_i$'s observed before time $j$, $j=1, \ldots, n$, whereas the probability of choosing a new tag depends only on the part of the stick that is left at time $n$. This is evident if we consider the alternative characterization of  \eqref{betagos1} with $p_{n,j}=W_j \prod_{i=j+1}^{n}(1-W_i)$ and $r_{n}(W_1,\dots,W_n)=\prod_{i=1}^n (1-W_i)$, $W_i \sim {\Be}(\beta_i, \alpha_i)$ and choose $\beta_i=1$ and $\alpha_i=\theta$ as in the DP.  Then, $p_{n,j}=W_j \prod_{i=j+1}^{n}(1-W_i)$, $j=1, \ldots, n$. For $n=3$, $p_{3,1}=W_1(1-W_2)(1-W_3)$, $p_{3,2}=W_2(1-W_3)$,  $p_{3,3}=W_3$. By contrast,  in a Dirichlet process each piece of the unitary stick is defined from what is left by the previous ones. \\

\section{A Beta-GOS hierarchical model}
\label{sec:gos-model}
In this Section, we show how the Beta-GOS process could be used as a
prior in a hierarchical model,  and we discuss a straightforward MCMC sampling algorithm for posterior inference.

\subsection{The hierarchical model.}
\label{sec:mode}
Beta-GOS processes can be used to model dependencies between non
exchangeable observations. Let
$\bm{Y}=(Y_1, \ldots, Y_m)^T$ be a vector of observations, e.g. a time
series. Then, following a Bayesian approach, we
can assume that the data can be described by a hierarchical model as
\begin{equation}\label{hiermodel1}
Y_i | \mu_i \stackrel{ind.}{\sim} p(y_i |\mu_i), \quad i=1, \ldots, m,
\end{equation}
for some probability density $p(\cdot | \mu_i)$, where the vector
$(\mu_1, \ldots, \mu_m)^T$ is a realization of a Beta-GOS process with
parameters $\alpha_i, \beta_i$ , $i=1, \ldots, m$, and base measure
$G_0$, which we succinctly denote as
\begin{align}\label{hiermodel2}
\begin{split}
\mu_1, \ldots, \mu_m &\sim \BGOS(\bm{\alpha}_m, \bm{\beta}_m, G_0),
\end{split}
\end{align}
i.e. is a sample from a random distribution characterized by the predictive rule \eqref{betagos1}, for some $W_i \sim \Be(\alpha_i, \beta_i)$, $i=1, \ldots, m$.
As noted in Section \ref{sec:gos-prior}, any Beta-GOS  defines a CID sequence.  In particular, marginally $\mu_i \sim G_0$, $i=1, \ldots, m$.
Therefore, $G_0$ can be regarded as a centering distribution,
as in DP mixture models: $G_0$ can  represent a vague parametric prior
assumption on the distribution of the parameters of interest.
The hierarchical model may be extended by putting hyper-priors on the remaining parameters of the model, including the parameters of the Beta-GOS $(\bm{\alpha}_m, \bm{\beta}_m, G_0)$, although here we focus on the characterization of the behavior of the Beta-GOS for fixed choices of the Beta parameters. 

We  conclude this Section by noting that the sequence $Y_1, Y_2, \ldots$, defined through \eqref{hiermodel1} and \eqref{hiermodel2}, with joint density
\[
\int \prod_{i=1}^m p(y_i |\mu_i) \pi(d\mu_1, \ldots, d\mu_m), \quad m \geq 0,
\]
and  $\pi(\cdot)\equiv \BGOS(\bm{\alpha}_m, \bm{\beta}_m, G_0)$, is also a CID sequence.
Therefore, although not exchangeable, the $Y_{n+j}$'s, $j \geq 1$ are conditionally identically distributed given $(Y_1, \ldots, Y_n, \mu_1, \ldots, \mu_n)$.
For a proof of this statement, see Proposition \ref{propCIDness} in the Appendix.\\

\subsection{MCMC posterior sampling.}
\label{sec:MCMCsampling}
Posterior inference  for the model
\eqref{hiermodel1}-\eqref{hiermodel2} entails learning about the
vector of random effects $\mu_i$ and their clustering structures. As
the posterior is not available in closed form, we need to revert to
MCMC sampling.
In this Section, we describe a Gibbs Sampler  that relies on
sampling the subsequent cluster assignments  of the observations $Y_1, \ldots, Y_m$ according to the rule \eqref{betagos1}.
To do this, the partition structure will be described by introducing a
sequence of labels $(C_n)_{n \geq 1}$ recording the pairing of each observation according to
\eqref{betagos1}, i.e. which other data point, among those with index
$j<i$,  the $i$th observation has been matched to. Hence, here
the label $C_i$ is not a simple indicator of the cluster  membership, as it is typical in most MCMC algorithms devised for the Dirichlet
process,  although cluster membership can be easily retrieved by
analyzing the sequence of pairings. In what follows, $C_i$ will be sometimes referred to as  the $i$-th pairing label.
In particular,  if the $i$-th observation is not paired to any
of those preceding, $C_i=i$; in this case, the $i$-th  point consists of a draw from  the
base distribution $G_0$, and thus generates a new cluster. This
slightly different representation of data points in terms of
data-pairing labels, instead of cluster-assignment labels, turns
useful to develop an MCMC sampling scheme for non-exchangeable
processes, as it has been thoroughly discussed in \cite{Blei09}, who
have shown that such representation allows for larger moves in the
state space of the posterior and faster mixing of the sampler.
It is easy to see that  the pairing sequence  $(C_n)_{n \geq 1}$
assigns $C_1=1$ and has distribution
\begin{equation}\label{eq:C_n}
\begin{split}
 P\{C_n=i|C_1,\dots,C_{n-1},W\}& =P\{C_n=i|W_1,\dots,W_{n-1}\} \\
&=r_{n-1} \J \{i=n \}+p_{n-1,i} \J\{i \not=n\},  \\
\end{split}
\end{equation}
for $i=1,\dots,n$,  where $\J(\cdot)$ denotes, as usual, the indicator function, such that, given a set A, $\J(i \in A)=1$ if $i \in A$ and $0$ otherwise.
As mentioned, the clustering configuration is a by-product of the
representation in terms of data-pairing labels. If two observations
are connected by a sequence of interim pairings, then they are in the
same cluster. Given $C=(C_1,\dots,C_m,\dots)$, let $\Pi(C)$ denote the
partition on $\N$ generated by $C$.  Accordingly, if $(\mu_k^*)_{k\geq
  1}$ is a sequence of  independent random variables with common
distribution $G_0$, we set $\mu_i=\mu^*_k$ if $i$ belongs to
$\Pi(C)_k$, i.e. the $k$-th block (cluster) of $\Pi(C)$.
For any $m$ and any $i\leq m$, let $C(m)=(C_1,\dots,C_m)$, $C_{-i}=(C_1,\dots,C_{i-1},C_{i+1},\dots,C_m)$; analogously, let $W(m)=(W_1,\dots, W_m)$, and  $W_{-i}=(W_1,\dots, W_{i-1},W_{i+1},\dots,W_m)$. Then, the full conditional for the pairing indicators $C_i$'s is
\begin{align}\label{eq:Gibbs1}
\begin{split}
P\{C_i=j|C_{-i},&Y(m),W(m)\}  \propto P\{C_i=j, Y(m)| C_{-i},W(m)  \} \\
& = P\{Y(m)| C_i=j, C_{-i},W(m)  \} P\{C_i=j| C_{-i},W(m)  \}. \\
\end{split}
\end{align}
The second term in \eqref{eq:Gibbs1} is the prior predictive rule \eqref{eq:C_n}, whereas
\[
 P\{Y(m)| C_i=j, C_{-i},W(m)  \} = \prod_{k=1}^{|\Pi{(C_{-i},j)}|} \int \prod_{ l \in \Pi{(C_{-i},j)}_k}
 p(Y_l|\mu_j^*) G_0(d\mu^*_j),
\]
where $\Pi{(C_{-i},j)}$ denotes the partition generated by
$(C_1,\dots,C_{i-1},j,C_{i+1},\dots,C_m)$. If $G_0$ and $p(y|\mu)$ are conjugate, the latter integral has a closed form solution. The non-conjugate case could be handled by appropriately adapting the algorithms of \cite{Mac98} and \cite{Neal2000}. Instead, we believe that split and merge moves as the ones considered in  \cite{Jain07} and \cite{Dahl05} are more problematic to implement given the implied exchangeability of the clustering assignments in those algorithms. As far as the full conditional for the latent variables $W_i$'s,  we can show that $W_i  |C(m),W_{-i},Y(m) \sim \text{Beta}(A_i, B_i)$, where $A_i=\a_i+\sum_{j= i+1}^m \J\{C_j<i \,\, \text{or}\,\, C_j=j  \}$, and $B_i=\beta_i+ \sum_{j=i+1}^m\J\{ C_j=i  \}$; hence, they depend on only on the clustering configurations and not on the values of $W_{-i}$.

Then, consider the set of cluster centroids $\mu_i^*$'s. The algorithm described so far allows faster mixing of the chain by integrating over the distribution of the $\mu_{i}^{*}$. However, in case we were interested on inference on the vector $(\mu_1, \ldots, \mu_m)$, it is possible to sample the unique values at each iteration of the Gibbs sampler, from
\begin{align}\label{xgivenpast}
\begin{split}
P\{\mu_j^*|C(m),W(m),Y(m)\}  \propto
\prod_{i \in \Pi_j(m)} p(Y_i|\mu^*_j)  G_0(d\mu^*_j),
\end{split}
\end{align}
where $\Pi_j(m)$ denotes the partition set of the observations such that $\mu_i=\mu_j^*$, $i=1, \ldots, m$. Again, if $p(y| \mu)$ and $G_0$ are conjugate, the full conditional of $\mu_j^*$ is available in closed form, otherwise we can update $\mu_{j}^{*}$ by standard Metropolis Hastings algorithms \citep{Neal2000}.

In addition, we note that if  a prior distribution for the Beta hyper-parameters $\bm{\alpha}_m$ and $\bm{\beta}_m$, say $\pi(\bm{\alpha}_m,\bm{\beta}_m)$, were to be specified, one could implement a Metropolis Hasting scheme to learn about their posterior distribution, since it can be shown that
\begin{align}\label{eq:mh4hyper}
P(\bm{\alpha}_m,\bm{\beta}_m| C(m), Y(m))\propto \pi(\bm{\alpha}_m,\bm{\beta}_m)\, \prod_{i=1}^m \frac{B(A_i, B_i)}{B(\alpha_i, \beta_i)},
\end{align}
where $A_i$ and $B_i$ are defined as above and $B(x, y)={\Gamma(x)\Gamma(y)}/{\Gamma(x+y)}$ denotes the Beta function. Equation \eqref{eq:mh4hyper} is an adaptation of well known results for the Dirichlet Process \citep{Esco95} to the Beta-GOS process. A thorough study of the efficiency of this algorithm, however, as well as the choice of adequate proposal distributions is beyond the scope of this work and will be pursued elsewhere.


\section{A Simulation study}
\label{sec:sim}

In this Section, we provide a full specification for model \eqref{hiermodel1}--\eqref{hiermodel2} and test the properties of the Beta-GOS prior on a set of simulated examples; more specifically, we develop some comparison with the Dirichlet Process and popular hidden Markov Models (HMM).

\subsection{Model specifications}

Throughout this Section, model \eqref{hiermodel1}--\eqref{hiermodel2} will be specified as follows. First,  we assume a Gaussian distribution for the observables, $Y_i \sim \text{Normal}(\mu_i,\tau^2)$. The base measure  $G_0$ is  also assumed to be normal, $\text{Normal}(\mu_0,\sigma_0^2)$, and  $\tau^2 \sim \IG(a_0,b_0)$.
 The parameters of the latent Beta reinforcements,  $W_i \sim \Be (\alpha_i,\beta_i)$, are separately indicated in each simulation and are chosen to allow for a range of prior beliefs on the clustering behavior of the process (see Section \ref{s:priorproperties}).
Details of the MCMC algorithm for posterior inference and parameter estimation in the Beta-GOS model are given in Appendix \ref{app:mcmc}.

\subsection{Model fitting and parameter estimation}
\label{sec:sim-res-1}

A first simulation study considers an ideal setting. We generate 1000 samples of 101 observations each from the Beta-GOS model \eqref{hiermodel1}--\eqref{hiermodel2}, with (a) $\a_n=n$, $\b_n=1$ and (b) $\a_n=3$, $\b_n=1$. The first 100 points are used for fitting purpose, whereas the $101st$ point is used to assess goodness of fit. Without much loss of generality, we fix $\mu_0=0$ and $\sigma_0=10$. We mimic the  typical scale observed in  the data analyzed in Section \ref{sec:app-aberrations} and set $\tau=0.25$ to distinguish the sample variability from the variability of the base measure. We fit the data  using  a Beta-GOS hierarchical model, with default Beta hyper-parameters $\a_n=n$, $\b_n=1$, and study how well we can recover the basic characteristics of the data under such specification. We assume $\tau^2 \sim \IG(2.004, 0.0063)$ in the model fitting. This choice of the Inverse-Gamma hyper-parameters allows $\tau^2$ to have mean around $0.25^2$ and relatively large variance.  In addition, we fit a  DP mixture model with concentration parameter $\theta=1$,  which on the basis of Proposition \ref{prop-Kn} (a) can be seen as compatible with the parameters used in our model. The mixture of DP  model is fit to data using the R package ``DPpackage'' \citep{Jara07}. In this framework,  the Dirichlet Process provides a convenient comparison; however, we should stress that, in general, the underlying exchangeability assumption may not  be appropriate to fully capture the dependency structure of the data generating process.
 
The results of this simple simulation study are summarized in Table \ref{tb:sim1}. Table  \ref{tb:sim1} reports summary statistics aimed at providing synthetic measures of the goodness of fit, namely the estimated number of clusters and the accuracy of cluster assignments, together with  a measure of predictive bias.  Following the  machine learning terminology for classification performance metrics,  we call accuracy the ratio of the correct cluster assignments with respect to the total of  assignments. We compute the predictive bias as follows: for each sample, and each MCMC output, we predict a new observation on the basis of the estimated parameters and the clustering configurations provided by the algorithm, say $Y_{pred}$. The prediction is  compared with  the original value, $Y_{101}$. The predictive bias is simply the average of $|Y_{101}-Y_{pred}|$, and can be regarded as a measure of how well the model can predict future observations.
Nearly all data points were assigned to the correct clusters. The Beta-GOS appears to compare favorably in terms of predictive bias, especially when the data incorporate a stronger dependency structure. Most of the error is intrinsic to the data
generating process. As typical of most Bayesian nonparametric
models, including the DP, the ability of the model and estimation
algorithms to recover the ground truth may be affected by the choice
of the relative magnitudes of the hyper-parameters $\sigma_0^2$ and
$\tau^2$.    The Supplemental Materials contain additional results for several specifications of the data generating mechanism as well as  several choices of the hyper-parameters for model fitting, confirming the above remarks. \ech

\subsection{Fitting mixture of Gaussians}
\label{sec:sim-res-2}

A second simulation study is designed to assess the robustness of the Beta-GOS framework to model mis-specifications: i.e., we fit the proposed non-exchangeable model to exchangeable data. The DP process provides a sensible baseline for this study. 
More in detail, we first generate 1,000 data sets (101 observations each) from a Normal mixture model with five components. The components' means are sampled from a $\text{Normal}(\mu_0=0, \sigma_0=10)$, whereas their standard deviation is set either to $\tau=0.25$ or $\tau=0.5$ to provide some insight into the robustness of the results to different levels of noise. The vector of  mixture components' weights is chosen as $\pi=	(0.2,0.35,0.15,0.1,0.2)^T$. We fit the data with a Dirichlet Process $(\theta=1)$, and a Beta-GOS process, with  a) $\alpha_n=\beta_n=1$, and  b) $\alpha_n=n$, $\beta_n=1$.
Case (a) corresponds to a process with short autocorrelation expected {\it a priori} and, asymptotically, a finite number of clusters, whereas case (b)  assumes that  the rescaled number of clusters,  $K_n/\log(n)$,  converges
to a $Gamma(1,1)$, and $E[K_n] \sim \log(n)$. The choice of hyper-parameters for the Inverse-Gamma on $\tau^2$ sets the mean around the true value and allows for a relatively large variance.
The results of the simulations are shown in Table \ref{tb:sim2}.  Overall, the Beta-GOS framework is quite robust to model mis-specifications. For the mixture of Gaussians data, accuracy of cluster assignments was high  (94\%), that is   better or comparable to that of the DP; correspondingly, parameters' estimates were close to the true parameter values.   
  For all processes, the accuracy decreases slightly with increasing level of  noise. In Figure \ref{fig:cluster_histograms}, we report the posterior distribution of the number of clusters for the three processes, for the case $\tau=0.25$. In accordance with the findings of Proposition \ref{prop-Kn}, we can see that if $\alpha_n=\beta_n=1$ the distribution is more concentrated around the mean and fewer clusters are generated in the fit. \\
Finally, we note that in our simulations, posterior inference for the Beta-GOS process seemed only minimally affected by the two different specifications of the parameters of the Beta weights. This consideration     confirms \ech the suggestion that using $\alpha_n=n+\theta-1$, $\beta_n=1$    represents \ech a default choice in many applications, where there is no {\it a priori} information to guide parameter choice. In this case, $\theta$ can be chosen or estimated similarly as what is routinely done for mixtures of DPs.     The Supplemental Materials contain additional results for several specifications of the model hyper-parameters, overall confirming the above remarks. \ech


\begin{table}[t!]
\caption{Summary statistics for the simulation study  in Section \ref{sec:sim-res-1}. The table compares the Beta-GOS and a Dirichlet Process model under different specifications of hyper-parameters when the data generating process is Beta-GOS.  }
{\footnotesize
\begin{center}
\scalebox{0.99}{ %
\begin{tabular}{|l|cc|cc|}
\hline 
\multicolumn{1}{|l|}{\textbf{Data Generating Process: }} & \multicolumn{2}{c|}{\textbf{Beta-GOS} } & \multicolumn{2}{c|}{\textbf{Beta-GOS} }\tabularnewline
 & \multicolumn{2}{c|}{$\alpha_{n}=n,\, \beta_{n}=1$} & \multicolumn{2}{c|}{$\alpha_{n}=3,\, \beta_{n}=1$}\tabularnewline
\hline 
\textbf{Model Fitting Method}  & \textbf{Beta-GOS }  & \textbf{Dir. Proc.}  & \textbf{Beta-GOS }  & \textbf{Dir. Proc.}\tabularnewline
 & \multicolumn{1}{c}{$\alpha_{n}=n$, $\beta_{n}=1$} & \multicolumn{1}{c|}{$\theta=1$} & \multicolumn{1}{c}{$\alpha_{n}=n$, $\beta_{n}=1$} & \multicolumn{1}{c|}{$\theta=1$}\tabularnewline
\hline 
\textbf{Number of Clusters}  &  &  &  & \tabularnewline
Ground Truth  & \multicolumn{2}{c|} {5.24$\pm$ 3.88}  
              & \multicolumn{2}{c|} {4.14 $\pm$ 1.81} \tabularnewline
Estimation  & 4.30$\pm$2.67  & 4.51$\pm$2.62  & 3.61 $\pm$ 1.49  & 3.96 $\pm$ 1.72 \tabularnewline
\hline 
Accuracy of Cluster Assignment  & 0.97$\pm$0.06  & 0.96$\pm$0.08  & 0.99 $\pm$ 0.01 & 0.99 $\pm$ 0.02 \tabularnewline
\hline 
Predictive Bias  & 4.13$\pm$7.18  & 4.34$\pm$7.27  & 0.67 $\pm$ 2.61 & 1.29 $\pm$ 3.93 \tabularnewline
\hline 
\end{tabular}} 
\end{center}

}
\label{tb:sim1}
\end{table}

\begin{sidewaystable}[htdp]
\caption{Summary statistics for the simulation study in Section \ref{sec:sim-res-2}. The table compares the Beta-GOS and a Dirichlet Process model under different specifications of hyper-parameters when the data generating process is a mixture of 5 gaussian components.}
{\footnotesize
\begin{center}
\scalebox{0.99}{
\begin{tabular}{|l|ccc|ccc|}
\hline 
 \multicolumn{7}{|l|}{\textbf{Data Generating Process: Gaussian Mixture} - 5 Gaussians}\tabularnewline
\hline 
\multicolumn{1}{|l|}{\footnotesize \bf True Sample Variability} & \multicolumn{3}{c|}{$\tau=0.25$} & \multicolumn{3}{c|}{$\tau=0.5$}\tabularnewline
\hline 
\textbf{Model fitting Method} & \multicolumn{2}{c}{\textbf{Beta-GOS }} & \textbf{Dir. Proc.} & \multicolumn{2}{c}{\textbf{Beta-GOS }} & \textbf{Dir. Proc.}\tabularnewline
 & $\alpha_{n}=n,\, \beta_{n}=1$ & $\alpha_{n}=\beta_{n}=1$ & $\theta=1$ & $\alpha_{n}=n,\, \beta_{n}=1$ & $\alpha_{n}=\beta_{n}=1$ & $\theta=1$\tabularnewline
\hline 
 Estimated Number of Clusters & 4.95$\pm$0.97 & 4.71$\pm$0.76 & 5.52$\pm$1.48 & 4.70$\pm$1.32 & 4.19$\pm$0.99 & 5.28$\pm$1.90\tabularnewline
\hline 
 Accuracy of Cluster Assignment & 0.94$\pm$0.09 & 0.93$\pm$0.09 & 0.93$\pm$0.09 & 0.86$\pm$0.11 & 0.84$\pm$0.12 & 0.85$\pm$0.13\tabularnewline
\hline 
 Predictive Bias & 8.86$\pm$9.02 & $8.73\pm$9.02 & 8.84$\pm$8.99 & 8.53$\pm$8.61 & 8.39$\pm8.37$ & 8.55$\pm$8.61\tabularnewline
\hline 
 Estimated Sample Variability & 0.25$\pm$0.01 & 0.25$\pm$0.01 & 0.27$\pm$0.05 & 0.56$\pm$0.68 & $0.69\pm1.19$ & 0.62$\pm$0.19\tabularnewline
\hline 
\end{tabular}}
\par\end{center}

}
\label{tb:sim2}
\end{sidewaystable}

\begin{sidewaystable}[htdp]
\caption{Summary statistics for the simulation studies described in Section \ref{sec:sim-res-3}. The table compares the Beta-GOS and a hidden Markov  model under different specifications of hyper-parameters. The data generating process assumes a hidden semi-Markov with state occupancy distribution  NegBin$(15,0.15)$ and two levels of the sampling noise $\tau=0.25$ and $\tau=0.5$.}
{\footnotesize
\begin{center}
\begin{tabular}{|lccc|ccc|}
\hline 
\multicolumn{7}{|l|}{\textbf{i) Data Generating Process: Hidden Semi Markov Model (HSMM)
with 4 states and }NegBin$(15,0.15)$}\tabularnewline
\hline 
\multicolumn{1}{|l|}{\textbf{Model Fitting Method}} & \multicolumn{3}{c|}{\textbf{Beta-GOS}} & \multicolumn{3}{c|}{\textbf{HMM}}\tabularnewline
 & \multicolumn{1}{|c}{$\alpha_{n}=n;\,\beta_{n}=1$} & \multicolumn{1}{c}{$\alpha_{n}=5;\,\beta_{n}=1$} & \multicolumn{1}{c|}{$\alpha_{n}=1;\beta_{n}=1$}  & \multicolumn{1}{c}{3 States} & \multicolumn{1}{c}{4 States} & \multicolumn{1}{c|}{5 States}\tabularnewline
\cline{1-7} 
 & \multicolumn{6}{c|}{$\tau=0.25$}\tabularnewline
 \hline 
Estimated Number of Clusters  & $3.89\pm 0.53$ & $4.04\pm0.59$ & $4.09\pm0.63$ & $2.98\pm0.13$ & $3.94\pm0.29$ & $4.87\pm0.55$ \tabularnewline
\hline 
Accuracy of Cluster Assignment  & $0.95\pm0.10$ & $0.97\pm0.07$ & $0.96\pm0.08$ & $0.72\pm0.09$ & $0.84\pm0.13$ & $0.90\pm0.13$\tabularnewline
\hline 
 & \multicolumn{6}{c|}{$\tau=0.5$}\tabularnewline
\hline 
Estimated Number of Clusters & $3.69\pm0.81$ & $3.89\pm0.96$ & $4.06\pm0.97$ & $2.99\pm0.12$ & $3.96\pm0.25$ & $4.90\pm0.48$ \tabularnewline
\hline 
Accuracy of Cluster Assignment & $0.86\pm0.14$ & $0.90\pm0.12$ & $0.90\pm0.12$ & $0.71\pm0.11$ & $0.83\pm0.12$ & $0.88\pm0.13$\tabularnewline
\hline 
\end{tabular}
\par\end{center}

\label{tb:sim3}
}
\end{sidewaystable}



\begin{figure}[tb]
\minipage{0.32\textwidth}
  \includegraphics[width=\linewidth]{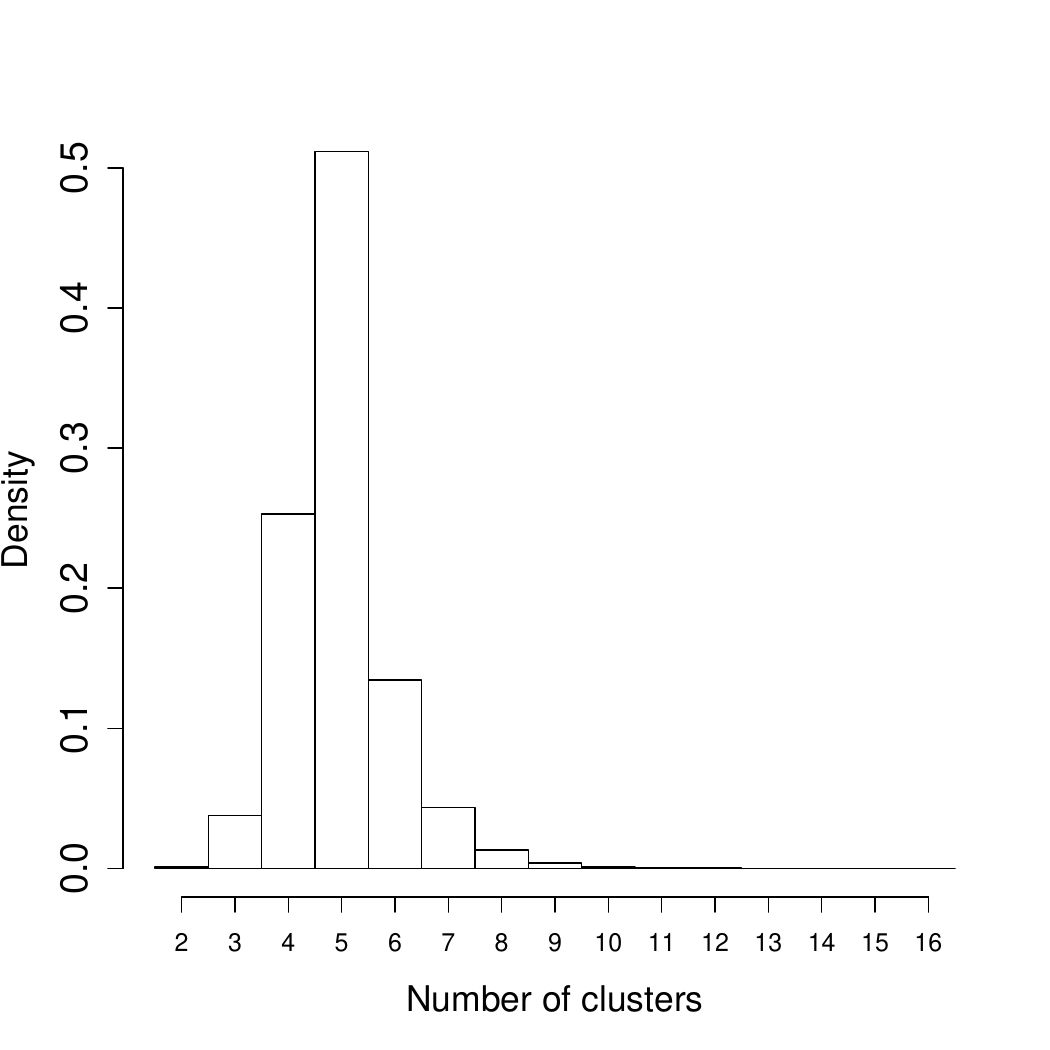}
	\subcaption[]{}
\endminipage\hfill
\minipage{0.32\textwidth}
  \includegraphics[width=\linewidth]{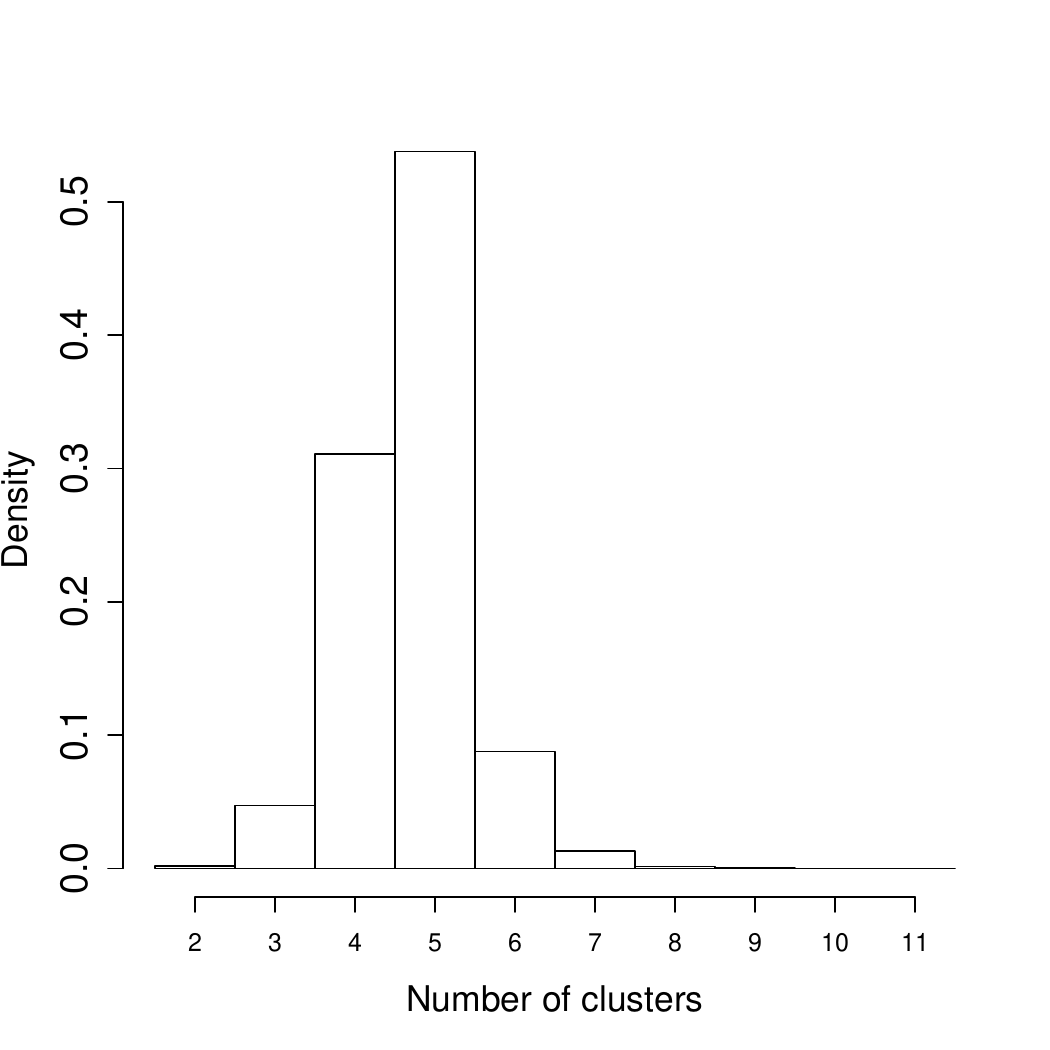}
  \subcaption[]{}
\endminipage\hfill
\minipage{0.32\textwidth}%
  \includegraphics[width=\linewidth]{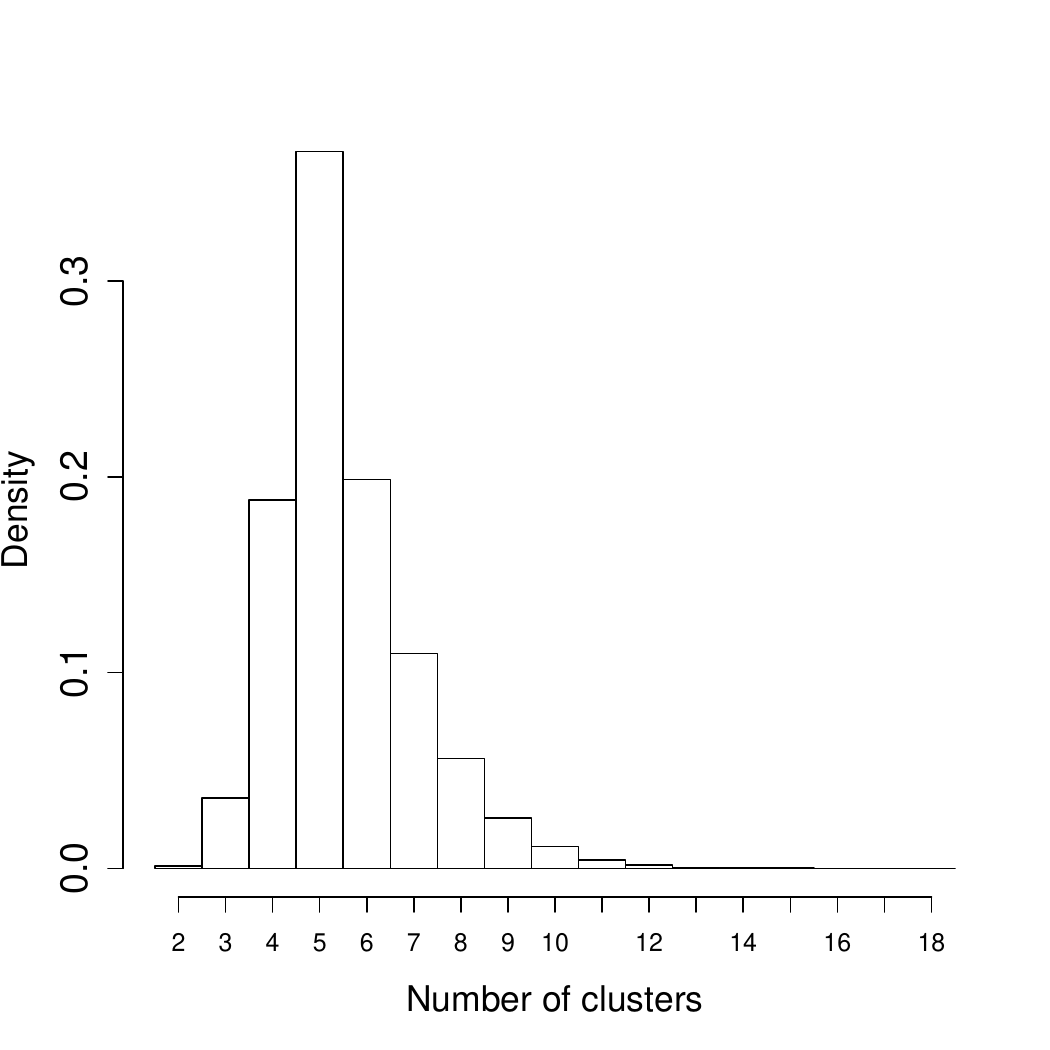}
  \subcaption[]{}
\endminipage
\caption{Posterior distribution of the number of clusters in the simulation of Section \ref{sec:sim-res-2} ($\tau=0.25$). Case (a) corresponds to a Beta-GOS($\alpha_n=n$, $\beta_n=1$), case (b) to a Beta-GOS($\alpha_n=\beta_n=1$)  and case (c) to a Dirichlet Process with parameter $\theta=1$.} \label{fig:cluster_histograms}
\end{figure}



\begin{figure}[tb]
\minipage{0.32\textwidth}
  \includegraphics[width=\linewidth]{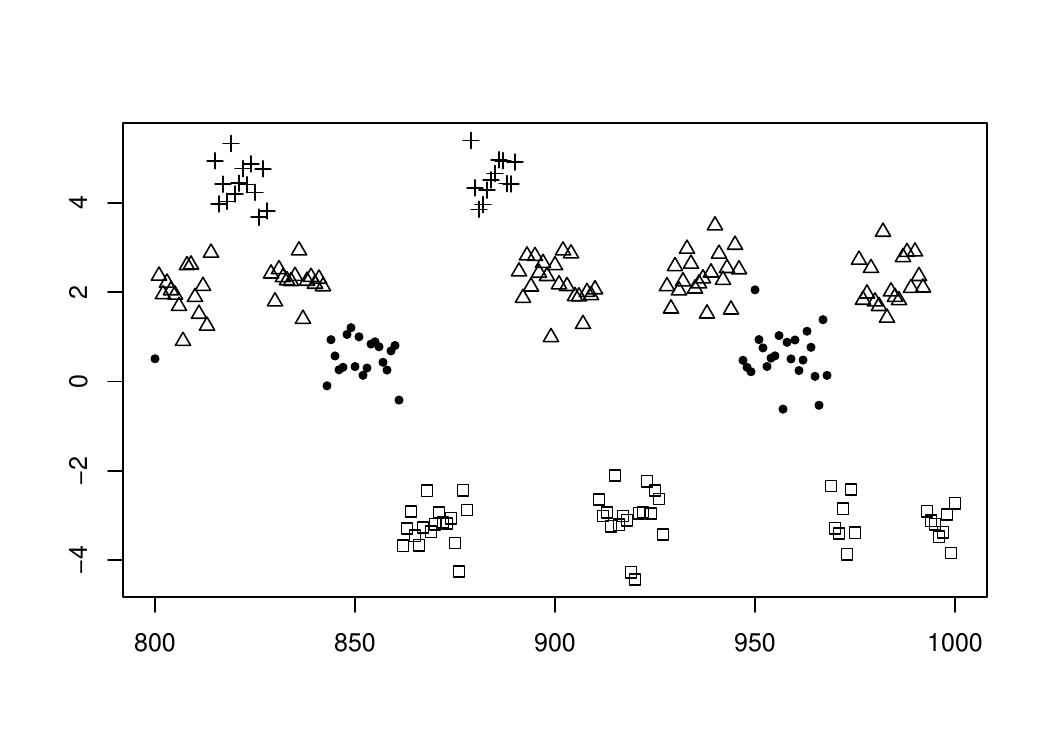}
\endminipage\hfill
\minipage{0.32\textwidth}
  \includegraphics[width=\linewidth]{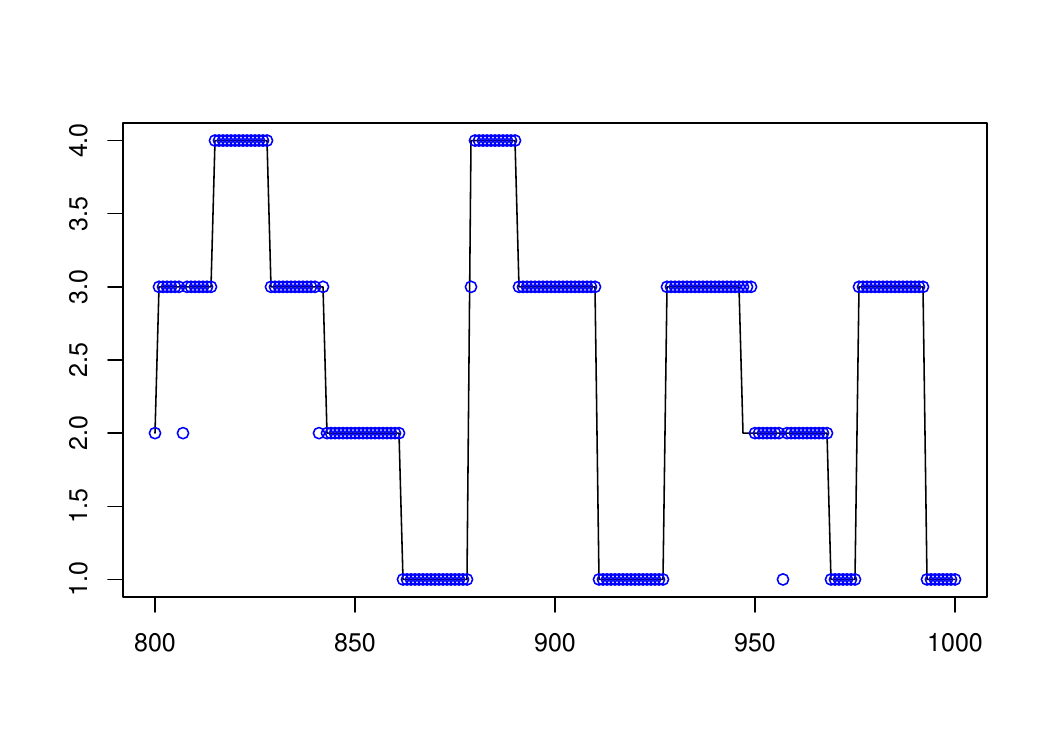}
\endminipage\hfill
\minipage{0.32\textwidth}%
  \includegraphics[width=\linewidth]{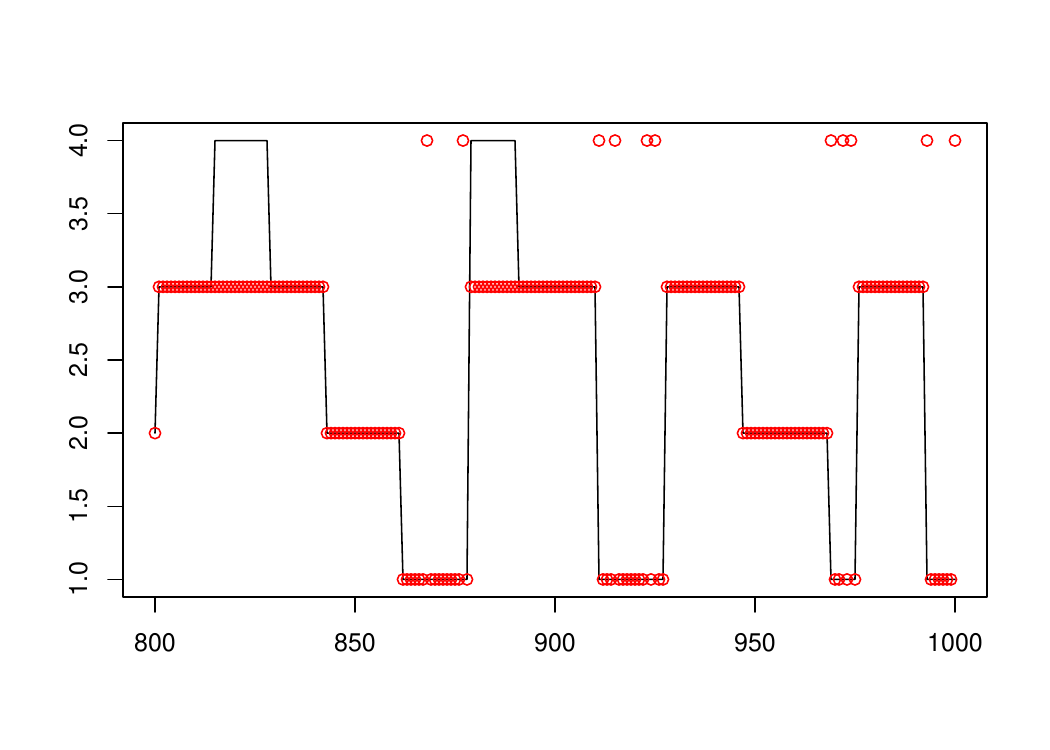}
\endminipage\\
\minipage{0.32\textwidth}
  \includegraphics[width=\linewidth]{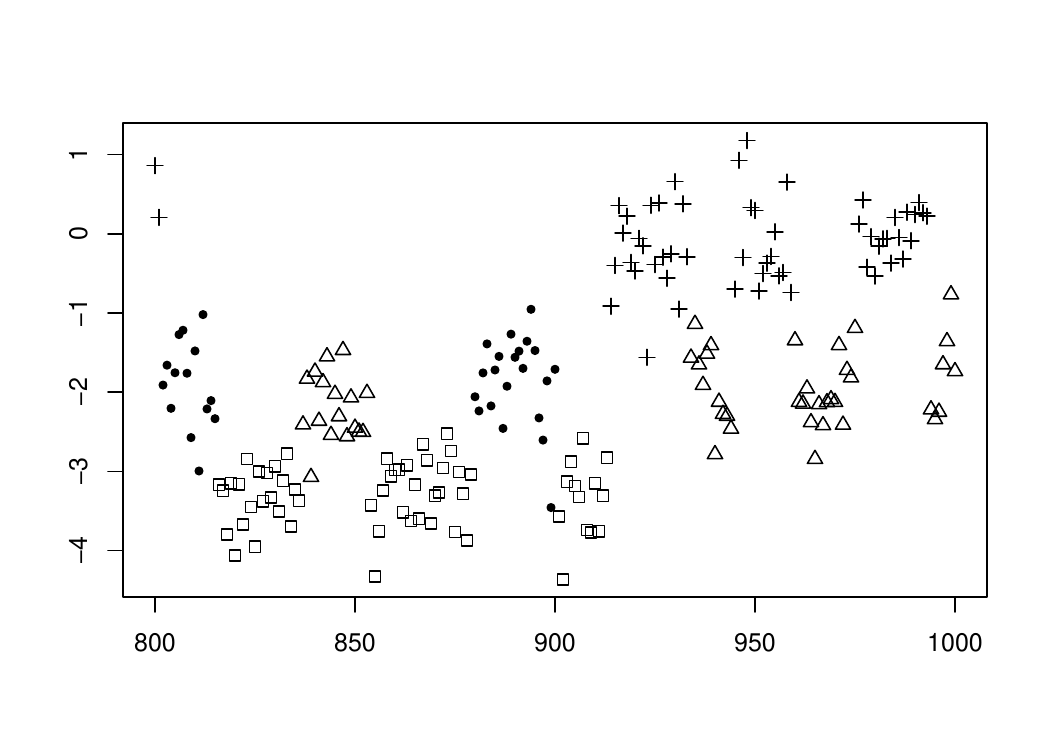}
	\subcaption[]{Ground Truth}
\endminipage\hfill
\minipage{0.32\textwidth}
  \includegraphics[width=\linewidth]{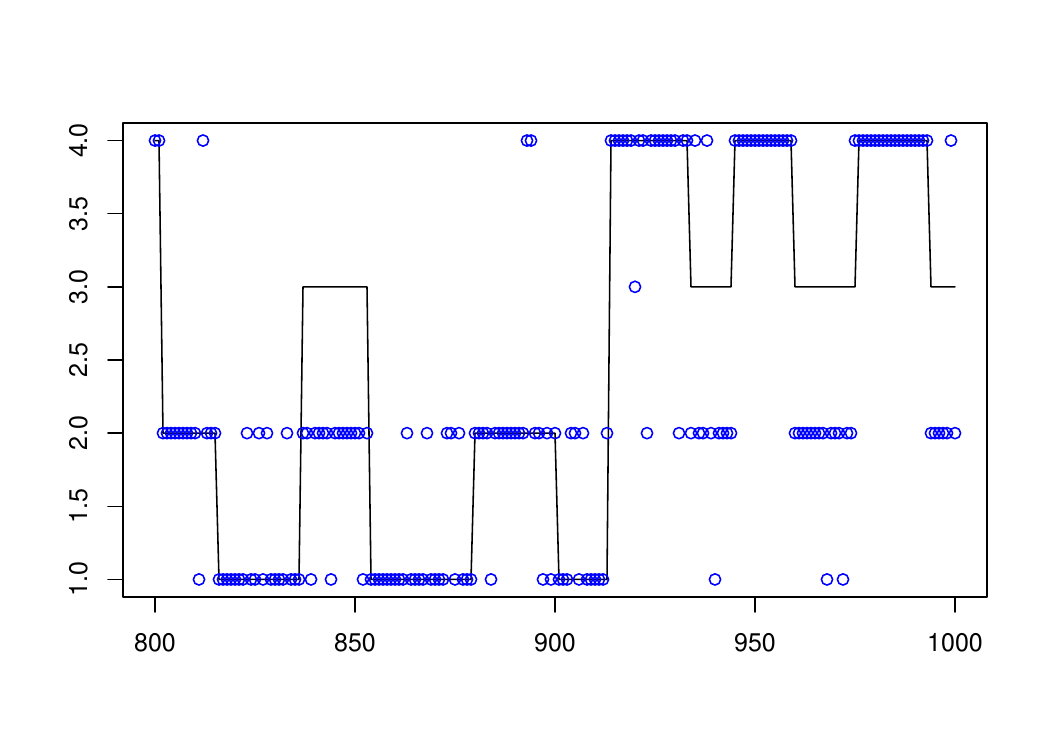}
  \subcaption[]{Beta Gos}
\endminipage\hfill
\minipage{0.32\textwidth}%
  \includegraphics[width=\linewidth]{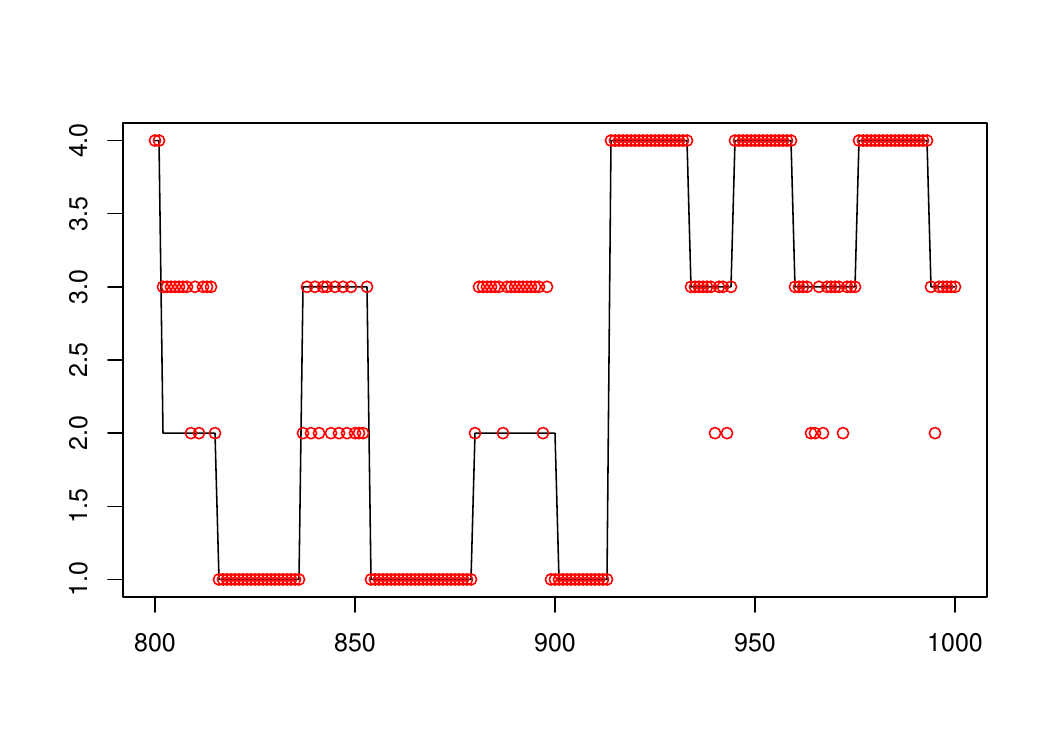}
  \subcaption[]{HMM}
\endminipage
\caption{   Illustrative segmentation-type plots for  the simulation study in Section 5.4. Column (a): subset of data for two replicates. Column (b) {\it top}: an example of allocation  for a Beta-Gos($\alpha_n=1,\beta_n=1)$ plotted {\it vs} the truth (black line); column (b) {\it bottom} considers a Beta-Gos($\alpha_n=n,\beta_n=1)$. Column (c) illustrates the fitting by a HMM with 4 states. \ech } \label{fig:Segmentation_plots}
\end{figure}

\subsection{Fitting Hidden Semi-Markov Models}
\label{sec:sim-res-3}
A third simulation study is designed to assess the robustness of the Beta-GOS framework to a mis-specification of a different nature. In many problems (e.g. change point detection), hidden Markov Models are used as computationally convenient substitutes for temporal processes that are known to be more complex than implied by first order Markovian dynamics. Here,  we generate non-exchangeable sequences from a hidden semi-Markov process \citep[HSMM;][]{Ferg80,Yu10} and study how the Beta-GOS process performs in fitting this type of data.  Hidden semi-Markov processes are an extension of the popular hidden Markov model where 
the time spent in each state (state occupancy or sojourn time) is  given by an explicit (discrete) distribution.  A geometric state occupancy distribution characterizes ordinary
hidden Markov models. Therefore, hidden semi-Markov process have  also been referred to as ``hidden Markov Models with explicit duration'' \citep{Mitchell95, Dewar12} or ``variable-duration hidden Markov Models'' \citep{Rabiner89}. 

We generate 1,000 datasets (1000 \ech observations each) using a hidden semi-Markov process with four states and a negative binomial distribution for the state occupancy distribution. More specifically, we parametrize the negative binomial in terms of its mean and an ancillary  parameter, which is directly related to the amount of overdispersion of the distribution \citep{Hilbe11, airoldi:2006}.   If the data are not overdispersed, the Negative Binomial reduces to the Poisson, and the ancillary parameter is zero.    For the simulations presented here, we consider a NegBin$(15,0.15)$,  which corresponds to assuming a large overdispersion (17.25). We also consider $\tau=0.25$ and $\tau=0.5$ in order to explore robustness to different levels of noise. 
We fit the data  by means of a Beta-GOS model with Beta hyper-parameters defined by:  a) $\alpha_n=n, \beta_n=1$;  b) $\alpha_n=5, \beta_n=1$; c)  $\alpha_n=1, \beta_n=1$. Based on Proposition \ref{prop-Kn}, those choices correspond to assuming  different clustering behaviors;  in particular, different expected number of clusters {\it a priori}. 
We then compare the Beta-GOS with the fit resulting from  hidden Markov models, assuming $3$, $4$ and $5$ states, respectively.
Results from the simulations are reported in Table \ref{tb:sim3}, where the HMM was implemented using the R package ``RHmm'' \citep{RHMMpackage}. Table \ref{tb:sim3} shows that the  Beta-GOS is a viable alternative to HMM, as it can provide  more accurate inference than a single hidden Markov model where the number of states is fixed a priori. As expected, higher levels of noise decrease the accuracy of the estimates, but the  reduction affects the fit of the Beta-GOS and hidden Markov Models similarly. Furthermore, the fit obtained with the Beta-GOS appears quite robust to the different choices of the hyper-parameters. Figure \ref{fig:Segmentation_plots} illustrates the clustering induced by the Beta-Gos and a 4-state HMM for a subset of the data generated in two specific simulation replicates. The middle column illustrates the allocation, respectively, from  a Beta-Gos($\alpha_n=1, \beta_n=1$) ({\it top}) and a Beta-Gos($\alpha_n=n, \beta_n=1$)({\it bottom}), whereas  column (c) illustrates the clustering attained by the HMM. Caution is necessary in order to avoid over-interpreting the results in the figure. Overall, the segmentation-plots suggest similarity in the allocations induced by the Beta-GOS and the HMM. In some instances, the Beta-GOS fit seems to  allow shorter stretches of contiguous identical states, as illustrated in the top row of Figure \ref{fig:Segmentation_plots}. On the other hand, when data are characterized by elevated intra-claster variability, as in the bottom  row of Figure \ref{fig:Segmentation_plots}, both the Beta-Gos and the HMM could fail to attain a fair representation of the true clustering structure of the data. Our practical experience suggests that the issue is more prominent for the ``default'' Beta-Gos($\alpha_n=n, \beta_n=1$) than for the ``informative''  Beta-Gos($\alpha_n=a, \beta_n=b$) formulations. This is in accordance with the discussion in Section  \ref{s:priorproperties} and, in particular, with the consideration that a Beta-Gos($\alpha_n=n, \beta_n=1$) should represent a long memory process where all previous observations are expected to contribute the same weight in \eqref{betagos1}. 
The Supplementary Materials contain results for a wider range of parameter settings, as well as different data generating mechanisms, confirming the results noted above.\ech

\section{Quantifying chromosomal aberrations in breast cancer}
\label{sec:app-aberrations}

We first apply the Beta-GOS to a classic dataset that has been used to link patterns of chromosomal aberrations to breast cancer pathophysiologies in the medical literature \citep{chin:2006}.
The raw data measure genome copy number gains and losses over 145 primary breast tumor samples, across the 23 chromosomes, obtained using BAC array Comparative Genomic Hybridization (CGH). More precisely, the measurements consist of $log_2$ intensity ratios obtained from the comparison of cancer and
normal female genomic DNA  labeled with distinct fluorescent dyes and co-hybridized on a microarray in the presence of Cot-1 DNA to suppress unspecific hybridization of repeat sequences \citep[see][]{Redon09}. 
The analysis of array CGH data presents some challenges, because data are typically very noisy and spatially correlated. More specifically, copy numbers gains or losses at a region are often associated to an increased probability of gains and losses at a neighboring region.
We use the  Beta-GOS model developed in the previous Sections to analyze and cluster clones with similar level of amplification/deletion, for each breast tumor sample and each chromosome in the dataset. 
For array CGH data, it is typical to distinguish regions with a normal amount of chromosomal material, from regions with single copy loss (deletion),  single copy gain and amplifications (multiple copy gains). Therefore, we present here the results of the analysis where the latent Beta hyper-parameters are set to $\alpha_n=3$  and  $\beta_n=1$, corresponding to $E(K_n)=4$ states for large n (see Section \ref{s:priorproperties}). We have also considered $\alpha_n=n$ and $\beta_n=1$, with no  remarkable differences in the results. We complete the specification of model \eqref{hiermodel1}--\eqref{hiermodel2} with a vague base distribution, $\text{Normal}(0, 10)$, and a vague inverse gamma distribution for $\tau$ centered around $\tau=0.1$. This choice of $\tau$ is motivated by the typical scale of the array CGH data and is in accordance with similar choices in the literature \cite[see, for example][]{Guha08}. 

Figure \ref{fig:chromosome8} exemplifies the fit to chromosome 8 on two tumor samples. The model is able to identify regions of reduced copy number variation and high amplification. Note how contiguous clones tend to be clustered together, in a pattern typical of these chromosomal aberrations. Figure \ref{fig:rawdata} replicates Figure 1 in \cite{chin:2006} and shows the frequencies of genome copy number gains and losses among all 145 samples plotted as a function of genome location. In order to identify a copy number aberration for this plot,  for each chromosome and sample, at each iteration we consider the cluster with lowest absolute mean and order the other clusters accordingly. The lowest absolute mean is chosen to identify the copy neutral state. Following \cite{Guha08} any other cluster is identified as a copy number gain or loss if its mean, say $\hat{\mu}_{(j)}$, is farther than a specified threshold from the minimum absolute mean, say $\hat{\mu}_{(1)}$, i.e. if $\hat{\mu}_{(j)}-\hat{\mu}_{(1)}>\epsilon$. We  experimented with a range of choices of $\epsilon$ in the range $[0.05,0.15]$ and used $\epsilon=0.1$ for the current analysis. Furthermore, if the mean of a cluster is above the mean of all declared gains plus two standard deviations, all genes in that cluster are considered  high level amplifications. We identify a clone with an aberration (or high level amplification) if it is such in more than 70\% of the MCMC iterations; then, we compute the frequency of aberrations and high level amplifications among all 145 samples, which are reported, respectively, at the top and bottom of Figure \ref{fig:rawdata}. As expected, the clusters identified by the model tend to be localized in space all over the genome. This feature may be facilitated by the increasingly low reinforcement of far away clones embedded in the Beta-GOS, and corresponds to the understanding that  clones that live at adjacent locations on a chromosome can be either amplified or deleted together due to the recombination process.

\begin{figure}[t!]
 \centering
  \includegraphics[width=\textwidth]{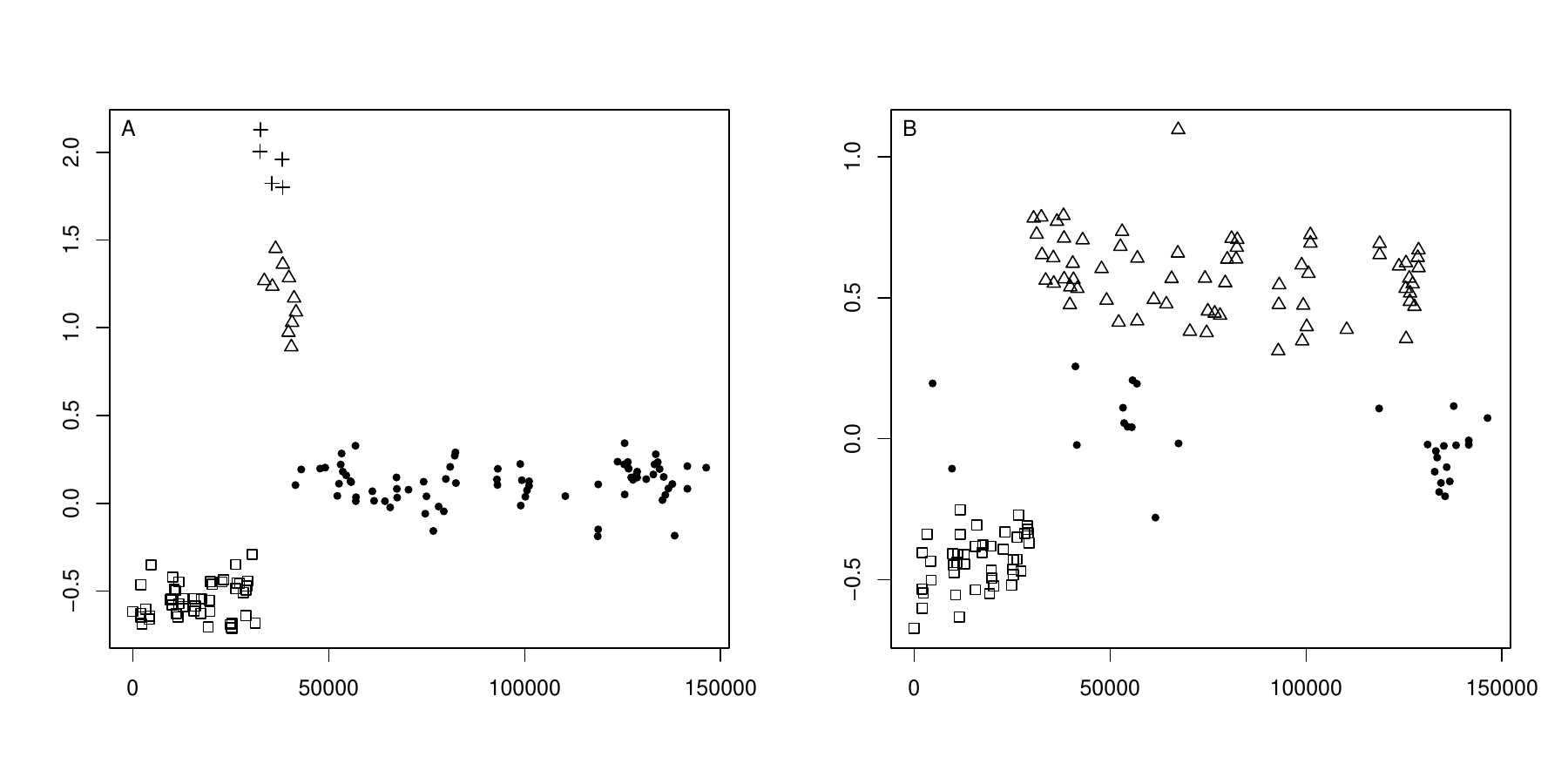}
\caption{Model fit overview:  Array CGH gains and losses on chromosome 8 for two samples of breast tumors in the dataset in \citep{chin:2006}. Points with different shapes denote different clusters.}
\label{fig:chromosome8}
\end{figure}

\begin{figure}[ht!]
 \centering
  \includegraphics[width=0.9\textwidth]{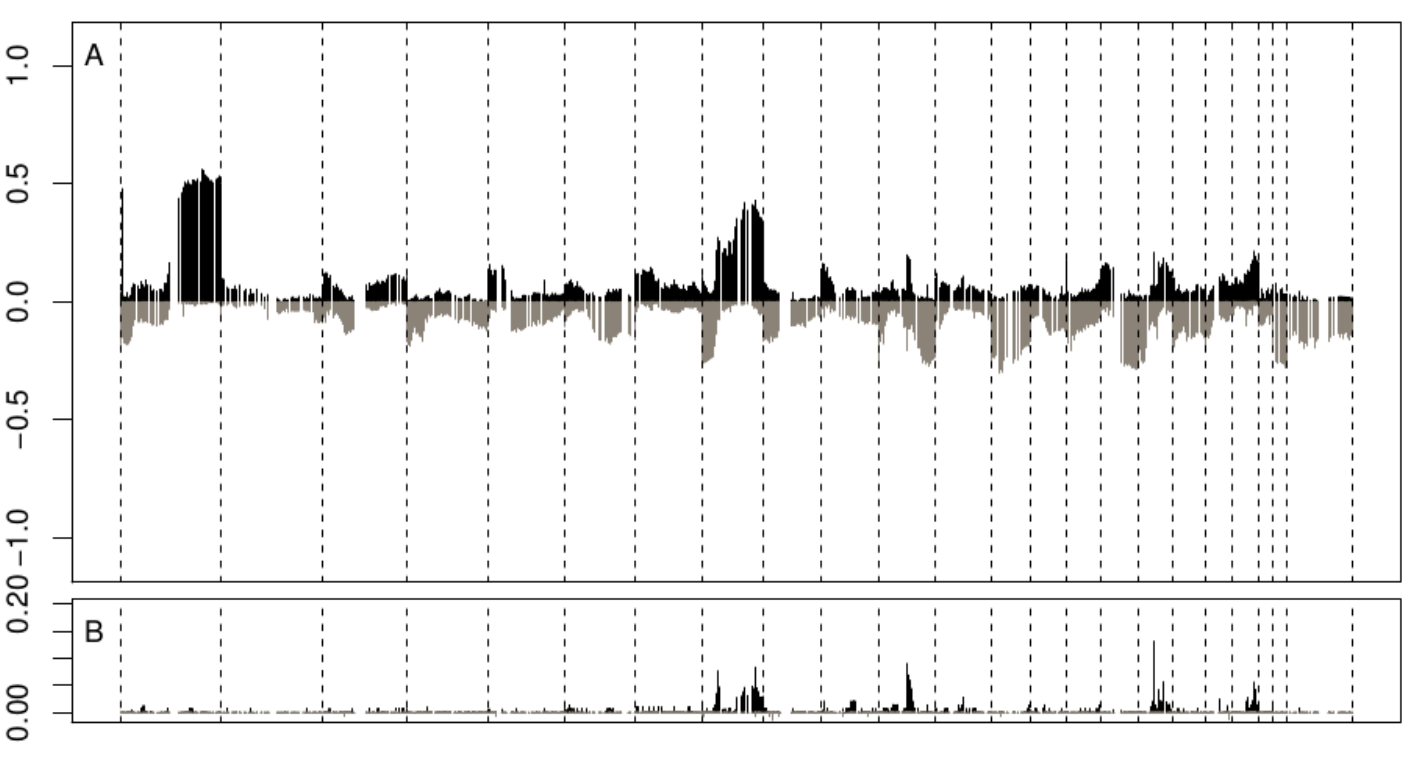}
\caption{A) Frequencies of genome copy number gains and losses plotted as a function of genomic location. B) Frequency of tumors showing high-level amplification. The dashed vertical lines separate the 23 chromosomes.}
\label{fig:rawdata}
\end{figure}

Finally, we considered some regions of chromosomes 8, 11, 17, and 20 that have been identified by \cite{chin:2006} and shown to correlate in their analysis to increased gene expression. We adapt the procedure described in \cite{Newton04} to compute a region-based measure of the false discovery rate (FDR) and determine the $q$-values for neutral-state and aberration regions estimated in our analysis. The $q$-value is the FDR analogue of the $p$-value, as it measures the minimum FDR threshold at which we may determine that a region corresponds to significant copy number gains or losses \citep{Sto03,Sto07}. More specifically, after conducting a clone based test as described in the previous paragraph, we identify regions of interest by taking into account the strings of consecutive calls. These regions then constitute the units of the subsequent cluster based FDR analysis. Alternatively, the regions of interest could be pre-specified on the basis of the information  available in the literature. The optimality of the type of procedures here described for cluster based FDR is discussed in \citealp{Sun12}. See also \citealp{Heller06}, \citealp{Muller07} and \citealp{Ji08}).
In Table \ref{tb:FDR} we report the $q$-values from a set of candidate oncogenes in well-known regions of recurrent amplification (notably, 8p12, 8q24, 11q13-14, 12q13-14, 17q21-24, and 20q13). Our findings confirm the previous detections of chromosomal aberrations in the same locations.

\begin{table}
\caption{False discovery rate analysis for clones with high-level amplification previously identified by \cite{chin:2006}. The individual amplicons are reported together with the locations of the flanking clones on the array platform. }
{\footnotesize
%

 \begin{center}
\begin{tabular}{| l | l | l |  l | l | l |     }
  \hline
{\bf Amplicon} & {\bf Flanking clone} & {\bf Flanking clone  }& \bf{ Kb }& {\bf Kb } & {\bf FDR} \\
   & {\bf (left)} & {\bf  (right) }& \bf{ start }& {\bf end} & {q-value}\\
\hline
   8p11-12 &RP11-258M15 & RP11-73M19& 33579 & 43001 & 0.021 \\ 
\hline
   8q24 & RP11-65D17&  RP11-94M13& 127186 & 132829& 0.021 \\ 
\hline
11q13-14 & CTD-2080I19& RP11-256P19 & 68482& 71659 &  0.022\\ 
\hline
  11q13-14 & RP11-102M18& RP11-215H8& 73337 & 78686 & 0.024\\ 
\hline
 12q13-14 &BAL12B2624 & RP11-92P22& 67191 & 74053 & 0.011\\ 
\hline
   17q11-12 & RP11-58O8& RP11-87N6& 34027& 38681 &  0.017\\ 
\hline
 17q21-24 & RP11-234J24&  RP11-84E24& 45775 & 70598 & 0.017\\ 
\hline
 20q13 & RMC20B4135&  RP11-278I13& 51669 & 53455& 0.021\\ 
\hline
  20q13 & GS-32I19& RP11-94A18& 55630 & 59444& 0.017\\ 
\hline

\end{tabular}
\end{center}

}
\label{tb:FDR}
\end{table}


Next, we apply our methodology to the analysis of a modern large-scale CGH array dataset \citep{Curtis12}. More specifically, here we consider one sample from the data published by \citet{Curtis12}. We  fit the Beta-GOS model to the entire sequence of 969,700 probes matched to genomics locations using a priority queue on the Harvard Odyssey cluster. Model fit took about 24 hours.
We also fit a Hidden Markov model and a Hidden semi-Markov model with Negative-Binomial run lengths times, both set to have three states \citep{Yau11}, to the same sample. The parameters of both models were estimated using standard techniques \citep{Rabiner89,Gued03}. The estimates for the Negative-Binomial run lengths, shared across states for simplicity, were $\hat r = 10$ and $\hat\pi\approx 0.25$, leading to a mean run length of 30 probes with a standard deviation of 10.

For validation purposes, we accessed a list of 152 consensus genomic locations where chromosomal aberrations were found in Breast cancer tumor samples. This list is included in the data files associated with The Genome Cancer Atlas (TCGA) project, partly curated by the Broad Institute and hosted by the NIH. The 152 consensus genomic locations range in size from 5 to 49 probes. 
This list provides a list of locations, which have been reported as likely altered in terms of DNA content in a number of publications, using multiple types of datasets and analyses.  Therefore, the list is  independent of the specifics of any particular technique, and it can be used as a reference for evaluating the comparative performance of our Beta-GOS model with Hidden Markov models, and Hidden semi-Markov models. 
For each method, we declared a success at detecting a chromosomal aberration (either deletion or amplification) at any of the 152 consensus genomic locations if the method correctly labeled at least 80\% of the probes associated with any given consensus genomic location. This choice was necessary since locations span multiple probes. According to this simple measure of performance, the Beta-GOS correctly labeled 133 locations (or 87.50\%) of the 152 consensus genomic locations as having a chromosomal aberration, versus 94 locations (or 61.84\%) using the Hidden Markov model, and 118 locations (or 77.63\%) using the Hidden semi-Markov model.
Of course, caution should be  should be taken against over-interpreting the results of a single illustrative example. However, the results from the simulation studies and  data analysis all concur to suggests that the Beta-GOS is a flexible model and that can be usefully employed in detecting chromosomal aberrations in array CGH data, since it can account for long range dependences in the sequence and achieve improved accuracy with respect to competing Hidden Markov Model based approaches.


\section{Concluding remarks.}
\label{sec:last}

Starting from the characterization of species sampling sequences in terms of their predictive probability functions, we have considered predictive rules where the weights are functions of latent Beta random variables.  The resulting  Beta-GOS process defines a novel and probabilistically coherent Bayesian Nonparametric model for non-exchangeable data.  We have discussed the clustering behavior of the Beta-GOS processes for some specifications of the latent Beta densities  and illustrated their use as priors in a hierarchical model setting. Finally, we have analyzed the performance of this modeling framework by means of a set of simulation studies. 
The results outlined in Section 6  illustrate how the
proposed Beta-GOS model can be a useful tool for the analysis of CGH array data.
In medical applications, for instance, it might be used to complement tumor sub-type
definition, or to suggest candidate genes for follow-up clinical studies.
We expect our approach will be useful in other applications where Hidden Markov
and semi-Markov model are currently considered as standard, e.g. in text segmentation and
speech processing \citep[e.g.,][]{ Rabiner89, BleiMoreno01, ChienFurui05, Ren10, YauHolmes13, Fox14}.\\
Recently, \cite{Teh06a},  \cite{Fox11}, and \cite{Yau11} have developed flexible and effective hierarchical Bayesian nonparametric extensions of hidden Markov models that allow posterior inference over the number of states. The Beta-GOS model provides an alternative, non-exchangeable, Bayesian nonparametric formalism  to model heterogeneity across non-exchangeable observations that are sequentially ordered, by enabling clustering in a number of unknown states.  Since the Beta-GOS model does not rely on the estimation of a single transition matrix across time points, as in the HMM,  we do not need to consider an explicit parameter to account for state persistence, as in \cite{Fox11}, or  assume a distribution for the sojourn times as in HSMMs. Indeed, since the predictive weights depend on the sequence of observations itself,  the Beta-GOS seems particularly convenient when the underlying generative process is non-stationary, e.g. as a possible alternative to more complicated non-homogeneous HMMs. \cite{Monteiro11} discuss a similar issue in a product partition model framework and explicitly assume that the observations in a  cluster have their distributions indexed by different parameters. Our approach is different, for example we do not need to explicitly model the dependence structure within the clusters. 

Arguably, the major obstacle we can foresee in the wider applicability of this type of models relies in the specification of the prior hyper-parameters in the latent Beta distributions.    Some specific suggestions have been provided in Section 3\ech. However,  in cases where there is not enough prior information to advise differently, our experience suggests that the default choice of the hyper-parameters outlined in Proposition \ref{prop-Kn}(a) not only reduces the problem to the choice of a single parameter as it is usual in DP mixture models, but  may also suffice for inferential purposes. Alternatively,  one could assume  a prior distribution on the  parameters of the Beta latent variables and conduct posterior inference by means of MCMC methods, as briefly discussed in Section \ref{sec:gos-model}.  Nevertheless, in specific applications the optimal modeling of the latent Beta densities requires further study and will be pursued elsewhere.  
 In addition, the proposed approach inherits the general computational limitations of nonparametric Bayesian methods. For example, a full MCMC algorithm for posterior inference may be unfeasible for genomic sequences with several millions of reads. Scalable algorithms may facilitate fast inference in those settings \citep[e.g.,][]{Colella07}. \ech \\
Finally, we believe that the flexibility of the latent specification and the possibility to tie the clustering implied by the Generalized P\'olya Urn scheme directly to a set of latent random variables gives an opportunity to further modeling the complex relationships typical of heterogenous datasets. For example, further  developments may substitute the general latent Beta specification with a probit/logistic specification, and define a Generalized P\'olya Urn scheme in the aims of \cite{Rodriguez10} that allows the clustering at each observation to be dependent on a set of (possibly sequentially recorded) covariates or curves. Similarly, we can imagine using multivariate Generalized P\'olya Urn schemes of the sort we describe in this paper to model time dependent parameters in time series, which may be important to identify time-varying structures or regime changes at the base of phenomena like the so called financial contagion, i.e. the co-movement of asset prices across global markets after large shocks \citep[see, for example,][]{Liu12}.

\section*{Acknowledgements}
\centering 
This work was partially supported by the National Science Foundation under grants DMS-1106980, IIS-1017967, and CAREER award IIS-1149662, by the National Institute of Health under grants R01 GM-096193 and P30-CA016672, and by the Army Research Office Multidisciplinary University Research Initiative under grant 58153-MA-MUR all to Harvard University. Fabrizio Leisen's research is supported by the European Community's Seventh Framework Programme FP7/2007-2013 under grant agreement number 630677. The authors would like to thank the Associate Editor and two anonymous referees  for suggestions that substantially improved the paper.



\begin{thebibliography}{100}
\newcommand{\enquote}[1]{``#1''}
\expandafter\ifx\csname natexlab\endcsname\relax\def\natexlab#1{#1}\fi


\bibitem[{Airoldi et~al.(2006)}]{airoldi:2006}
Airoldi E.~M., Anderson, A., Fienberg, S.E., Skinner, K.K. (2006)
\newblock Who wrote Ronald Reagan's radio addresses?
\newblock \emph{Bayesian Analysis}, 1, 289--320.
\bibitem[{Aoki, M.(2008)}]{Aoki08} Aoki M. (2008) Thermodynamic limits of macroeconomic or financial models: One- and two-parameter Poisson-Dirichlet models, \emph{Journal of Economic Dynamics and Control}, Elsevier, vol. 32(1), pages 66--84.
\bibitem[{Baladandayuthapani et~al.(2010)}]{Baladandayuthapani10} Baladandayuthapani V., Ji Y. Nieto-Barajas, L.E. and Morris, J.S. (2010) Bayesian random segmentation models to identify shared copy number aberrations for array CGH data. \emph{Journal of the American Statistical Association}, 105, 1358--1375.
\bibitem[{Bassetti et~al.(2010)}]{bas-cri-lei} Bassetti
 F., Crimaldi I. and Leisen F. (2010) {Conditionally identically
    distributed species sampling sequences.}
\emph{ Adv. in Appl. Probab.} 42, 433-459.
 \bibitem[{Berti et~al.(2004)}]{be-pra-ri} Berti P.,
   Pratelli L. and Rigo P. (2004) { Limit Theorems for a Class of
  Identically Di\-stri\-bu\-ted Random Variables.} \textit{ Ann.
  Probab.} {\bf 32} 2029--2052.
\bibitem[{Blackwell and MacQueen(1973)}]{black-mac} Blackwell D. and MacQueen J.B. (1973) { Ferguson distributions via P\'olya urn
    schemes.}   \textit{Ann. Statist.} {\bf 1}, 353--355.
\bibitem[{Blei and Frazier(2011)}]{Blei09} Blei D. and Frazier P. (2011) Distance dependent Chinese restaurant processes. \emph{Journal of Machine Learning Reseach}, 12:2461--2488. 
\bibitem[{Blei and Moreno(2001)}]{BleiMoreno01} Blei D.M., Moreno P. (2001) {Topic segmentation with an aspect Hidden Markov Model}, \textit{Proceedings of the 24th annual international ACM SIGIR conference}, 343--348.
\bibitem[{Cardin et~al.(2011)}]{Cardin11} Cardin N., Holmes, C., The Wellcome Trust Case Control Consortium, Donnelly P. and Marchini J. (2011) Bayesian Hierarchical Mixture Modeling to Assign Copy Number from a Targeted CNV Array. \emph{Genetic Epidemiology}, 35, 536--548.
\bibitem[{Charalambides(2005)}]{Charalambides} Charalambides C. A. (2005) Combinatorial methods in discrete distributions. Wiley Series in Probability
and Statistics. Wiley-Interscience, John Wiley \& Sons, Hoboken, NJ.
\bibitem[{Chin et~al.(2006)}]{chin:2006} Chin K., DeVries S., Fridlyand J.,  Spellman P.~T., Roydasgupta R., Kuo W.~L.,
  Lapuk A., Neve R.~M., Qian  Z., Ryder  T., Chen F., Feiler H., Tokuyasu T.,
  Kingsley C., Dairkee S., Meng Z., Chew K.,Pinkel  D., Jain A., Ljung B.~M.,
  Esserman L.,  Albertson D.~G.,  Waldman F.~M., and Gray J.~W. (2006)
\newblock Genomic and transcriptional aberrations linked to breast cancer
  pathophysiologies.
\newblock \emph{Cancer Cell}, 10\penalty0 (6):\penalty0 529--541.
\bibitem[{Chien and Furui(2005)}]{ChienFurui05} Chien J.T., Furui S. (2005) {Predictive hidden Markov model selection for speech recognition}, \textit{IEEE Transactions on Speech and Audio Processing}, {\bf 13}, 377--387.
\bibitem[{Colella et.~al.(2007)}]{Colella07} Colella S., Yau C., Taylor J.M., Mirza G., Butler H., Clouston P., Basset A.S., Seller A., Holmes C.C. and Ragoussis J. (2007) {QuantiSNP: an Object Bayes Hidden-Markov Model to detect and accurately map copy number variation using SNP genotyping data}   \textit{Nucleic Acids Research} {\bf 35}, 2013--2025.
\bibitem[{Curtis et~al.(2012)}]{Curtis12} Curtis C., Shah S.P., Chin S.-F., Turashvili G., Rueda O.M., Dunning M.J., Speed D., Lynch A.G., Samarajiwa S., Yuan Y., Gr\"af S., Ha G., Haffari G., Bashashati A., Russell R., McKinney S., METABRIC Group, Langer¿d A., Green A., Provenzano E., Wishart G., Pinder S., Watson P., Markowetz F., Murphy L., Ellis I., Purushotham A., B¿rresen-Dale A.-L., Brenton J.D., Tavar'e S., Caldas C., Aparicio S. (2012) {The genomic and transcriptomic architecture of 2,000 breast tumours re-veals novel subgroups}, \textit{Nature}, {\bf 486}, 7403, 346--352.
\bibitem[{Dahl(2005)}]{Dahl05} Dahl D. B. (2005) Sequentially-Allocated Merge-Split Sampler for Conjugate and Nonconjugate Dirichlet Process Mixture Models, \emph{Technical Report.}
\bibitem[{Dahl et al.(2008)}]{Dhal08} Dahl D. B., Day  R. and  Tsai J. W. (2008) {Distance-Based Probability Distribution on Set Partitions with Applications to Protein Structure Prediction}, \emph{Technical Report}.
\bibitem[{DeSantis et~al.(2009)}]{DeSantis09} DeSantis S. M., Houseman E. A., Coull B.A., Louis D.N., Mohapatra G., Betensky, R.A. (2009) A Latent Class Model with Hidden Markov Dependence for Array CGH Data, \emph{Biometrics}, {\bf 65}, 4, 1296--1305.
\bibitem[{Dewar et~al.(2012)}]{Dewar12} Dewar M., Wiggins C., Wood F. Inference in Hidden Markov Models with Explicit State Duration Distributions, \emph{Signal Processing Letters, IEEE}, 19, 4, 235--238.  
\bibitem[{Du et~al(2010)}]{Du10}  Du L.,  Chen M.,   Lucas, J. and   Carin, L. (2010)  Sticky Hidden Markov Modeling of Comparative Genomic Hybridization, \emph{IEEE Transactions on Signal Processing}, {\bf 58}, 10, 5353--5368.
\bibitem[{Escobar and West (1995)}]{Esco95} Escobar M. and West M. (1995) Bayesian Density Estimation and Inference Using Mixtures, \emph{Journal of the American Statistical Association}, {\bf 90}, 577-588.
\bibitem[{Ferguson(1980)}] {Ferg80} Ferguson J. D. (1980), Variable duration models for speech, \emph{Proceedings of
the Symposium on the Applications of Hidden Markov Models to Text and
Speech}, 143--179.
\bibitem[{Fortini et.~al(2000)}]{Fortini2000} Fortini S., Ladelli  L. and  Regazzini E. (2000)  Exchangeability,
predictive distributions and parametric models.  \emph{Sankhya}
Ser. A, {\bf 62}, no. 1, 86--109.
\bibitem[{Fox et al.(2011)}]{Fox11} Fox E.B., Sudderth E.B., Jordan M.I. and Willsky A.S. (2011) A sticky HDP-HMM with application to speaker diarization, \emph{Annals of Applied Statistics}, {\bf 5}, 2A, 1020--1056. 
\bibitem[{Fox et~al.(2014)}]{Fox14} Fox E.B., Hughes M.C., Sudderth E.B., Jordan M.I. (2014) {Joint modeling of multiple time series via the Beta process with application to motion capture segmentation}, \textit{Annals of Applied Statistics}. To appear.
\bibitem[{Griffiths(2007)}] {Griffiths07} Griffiths T.L., Sanborn   A.N., Canini K.R.,  and Navarro  D.J. (2007) Categorization as nonparametric Bayesian density estimation, M. Oaksford and N. Chater (Eds.), \emph{The Probabilistic Mind: Prospects for Rational Models of Cognition}, Oxford: Oxford University Press.
\bibitem[{Guedon(2003)}]{Gued03} Guedon Y. (2011) {Estimating hidden semi-Markov chains from discrete sequences}, \textit{Journal of Computational and Graphical Statistics}, {\bf 12}, 3, 604--639.
\bibitem[{Guha et~al.(2008)}]{Guha08} Guha S., Li Y. and Neuberg D. (2008). Bayesian Hidden Markov Modeling of Array CGH Data. \emph{Journal of the American Statistical Association}, 103, 485--497.
\bibitem[{Guha(2010)}]{Guha10} Guha S. (2010) Posterior Simulation in Countable Mixture Models for Large Datasets. \emph{Journal of the American Statistical Association},  {\bf 105}, 490, 775--786.
\bibitem[{Hansen and Pitman(2000)}]{han-pitman} Hansen B. and Pitman J.
  (2000) {Prediction rules for exchangeable sequences related to
    species sampling.}  \textit{Statist. Probab. Lett.} {\bf 46}
  251--256.
 \bibitem[{Heller et al.(2006)}]{Heller06} Heller  R., Stanley D.,  Yekutieli D.,  Rubin N., and  Benjamini Y. (2006). Cluster-based analysis of fmri data. \emph{Neuroimage}, 33, 599--608.
\bibitem[{Hjort et al.(2010)}]{Hjort10} Hjort  N.L., Holmes C., M\"uller  P. and Walker S.G. (2010) \emph{Bayesian Nonparametrics},  Cambridge University Press.
\bibitem[{Hilbe(2011)}]{Hilbe11} Hilbe J.M. (2011) \emph{Negative Binomial Regression},  Cambridge University Press.
\bibitem[{Ishwaran and Zarepour(2003)}]{IshZar03} Ishwaran H. and Zarepour M. (2003) {Random probability measures via P\'olya sequences: revisiting the Blackwell-MacQueen urn scheme}.  \verb+http://arxiv.org/abs/math/0309041+
\bibitem[{Jain and Neal(2007)}]{Jain07} Jain S. and Neal R.M. (2007)
Splitting and Merging Components of a Nonconjugate Dirichlet Process Mixture Model, \emph{Bayesian Analysis}, 3, 445--472.
\bibitem[Jara A. (2007)]{Jara07} Jara A. (2007) Applied Bayesian Non- and Semi-parametric Inference using DPpackage,  Rnews, {\bf 7}, 3, 17--26.
\bibitem[{Jbabdi et~al.(2009)}] {Jbabdi09} Jbabdi S., Woolrich M.W. and Behrens T.E.J. (2009) Multiple-subjects connectivity-based parcellation using hierarchical Dirichlet process mixture models, \emph{NeuroImage}, {\bf 44}, 2, 373--384.
\bibitem[Ji et~al.(2008)]{Ji08} Ji Y., Lu Y and Mills G. (2008) Bayesian models based on test statistics for multiple hypothesis testing problems. \emph{Bioinformatics}, {\bf 24}(7), 943-939.
\bibitem[{Kim et~al.(2006)}]{Kim06} Kim S., Tadesse M.G. and Vannucci M. (2006) Variable selection in clustering via Dirichlet process mixture models. \emph{Biometrika}, {\bf 93}(4), 877--893.
\bibitem[{Kingman(1978)}]{King78}  Kingman J. F. C. (1978) The representation of partition structures. \emph{J. London Math. Soc.} (2), 18(2), 374-380.
\bibitem[{Lee et al.(2008)}]{Lee08} Lee J., Quintana F., M\"uller P. and Trippa L. (2013) Defining Predictive Probability Functions for Species Sampling Models.  \emph{Statist.Sci.}, {\bf 28}, 2, 209--222.
\bibitem[{Liu et. al(2012)}]{Liu12} Liu Z., Windle J. and Scott J.C. (2012) The partition problem: case studies in Bayesian
screening for time-varying model structure. Technical report. Currently available at \verb+ http://arxiv.org/pdf/1111.0617.pdf+
\bibitem[{MacEachern and M\"uller(1998)}]{Mac98} MacEachern S. N. and M\"uller P. (1998) Estimating Mixture of Dirichlet Process Models, \emph{Journal of Computational and Graphical Statistics}, 7, 223--238.
\bibitem[{Marioni et al.(2006)}]{Marioni06} Marioni J.C., Thorne N.P., Tavare S., Radvanyi F. (2006) BioHMM: A heterogeneous hidden Markov model for segmenting array CGH data. \emph{Bioinformatics},22, 1144--1146.
\bibitem[{Mitchell et al.(1995)}]{Mitchell95} Mitchell C., Harper M., Jamieson L. , On the complexity of explicit duration HMMs, \emph{IEEE Transactions on Speech and Audio Processing}, 3, 2, 213--217.
\bibitem[Monteiro et al.(2011)]{Monteiro11} Monteiro, J.V., Assun\c{c}ao and Loschi, R.H. (2011) 
Product partition models with correlated parameter, \emph{Bayesian Analyis}, {\bf 6}, 4, 691--726.
\bibitem[{M\"uller et~al.(2007)M\"uller, Parmigiani and Rice}]{Muller07}
M\"uller P., Parmigiani G. and Rice K. (2007) FDR and Bayesian multiple
  comparisons rules.
\newblock In \textit{Bayesian Statistics 8} (eds. J.~Bernardo, M.~Bayarri,
  J.~Berger, A.~Dawid, Heckerman, A.~D., Smith and M.~West). Oxford, UK: Oxford
  University Press.
\bibitem[{M\"uller and Quintana(2010)}]{Mueller10} M\"uller P. and Quintana F. (2010) Random partition models with regression on covariates, \emph{Journal of Statistical Planning and Inference}, {\bf 140}, 10, 2801--2808.
\bibitem[{Navarro et al.(2006)}] {Navarro06} Navarro D.J., Griffiths T.L., Steyvers M. and Lee M.D. (2006) {Modeling individual differences using Dirichlet processes}. \emph{Journal of Mathematical Psychology}. In Special Issue on Model Selection: Theoretical Developments and Applications, Vol. 50, No. 2., pp. 101--122.
\bibitem[{Neal(2000)}]{Neal2000} Neal R.M. (2000) Markov Chain Sampling Methods for Dirichlet Process Mixture Models\emph{Journal of Computational and Graphical Statistics}, 9, 249--265.
\bibitem[{Newton et al.(2004)}]{Newton04} Newton M. A., Noueiry A., Sarkar D.  and Ahlquist P. (2004) Detecting differential gene expression with a semiparametric hierarchical mixture method. \emph{Biostatistics}, {\bf 5}, 2, 155---176.
\bibitem[{Park and Dunson(2007)}]{ParkDunson07} Park J.H. and Dunson D.B. (2010)
Bayesian generalized product partition model. \emph{Statistica Sinica}, {\bf 20}, 1203--1226
\bibitem[{Pitman(1996b)}]{Pit96b} Pitman J. (1996) Some developments of the Blackwell-MacQueen urn scheme. In T.S. Fer-guson et al., editor, \emph{Statistics, Probability and Game Theory; Papers in honor of
David Blackwell}, volume 30 of Lecture Notes-Monograph Series, pages 245-267. Institute of Mathematical Statistics, Hayward, California.
\bibitem[{Pitman(2006)}]{Pit06} Pitman J. (2006) \emph{Combinatorial Stochastic Processes.} Ecole  d'Et\'e
Probabilit\'es  de Saint-Flour XXXII  2002, Lecture Notes in Mathematics, Springer:Berlin / Heidelberg.
\bibitem[{Rabiner(1989)}]{Rabiner89} Rabiner L. R.  (1989) A Tutorial on Hidden Markov Models and Selected Applications in Speech Recognition, \emph{Proceedings of the IEEE}, 77, 2, 257--287. 
\bibitem[{Redon et al.(2009)}]{Redon09} Redon R., Fitzgeral T. and Carter, N.P. (2009) Comparative Genomic Hybridization: DNA labeling, hybridization and detection.
\emph{Methods Mol Biol.}  529: 267--278. 
\bibitem[{Ren et~al.(2010)}]{Ren10} Ren L., Dunson D., Lindroth S., Carin L. (2010) {Dynamic nonparametric Bayesian models for analysis of music}, \textit{Journal of the American Statistical Association}, {\bf 105}, 490, 458--472.
\bibitem[{Rodriguez et~al.(2010)}]{Rodriguez10} Rodriguez A., Dunson D.B. (2011) Nonparametric Bayesian models through probit stick-breaking processes. \emph{Bayesian Analysis}, {\bf 6}, 1,  145--178.
\bibitem[{Storey(2003)}]{Sto03}
Storey J. (2003) The positive false discovery rate: a Bayesian Interpretation and the q-value
\newblock \textit{The Annals of Statistics}, \textbf{31}, 6, 2013--2035.
\bibitem[{Storey et~al.(2007)Storey, Dai and Leek}]{Sto07}
Storey J., Dai J. and Leek J. (2007) The optimal discovery procedure for
  large-scale significance testing, with applications to comparative microarray
  experiments.
\newblock \textit{Biostatistics}, \textbf{8}, 414--432.
\bibitem[{Sudderth and Jordan(2009)}] {Sudderth09} Sudderth E. B. and Jordan M. I. (2009) Shared segmentation of natural scenes using dependent {P}itman-{Y}or processes. In \emph{ Neural Information Processing Systems} 22.
\bibitem[{Sun et. al(2014)}]{Sun12} Sun W., Reich B., Cai T., Guindani M., Schwartzman A. (2014) False Discovery Control in Large Scale Spatial Multiple Testing. \emph{Journal of the Royal Statistical Society (Series B).} To appear.
\bibitem[{Taramasco and Bauer(2012)}]{RHMMpackage} Taramasco O. and Bauer S. (2012) RHmm: Hidden Markov Models simulations and estimations.
\bibitem[{Teh et~al.(2006a)}] {Teh06a} Teh Y. W.,  Jordan M. I., Beal M. J. and Blei D. M. (2006a)
{Hierarchical {D}irichlet processes.}, \emph{J. Amer. Statist. Assoc.}, {\bf 101}, no. 476, 1566--1581.
\bibitem[{Teh(2006b)}]{Teh06b} Teh Y. W. (2006b) {A {H}ierarchical {B}ayesian {L}anguage model based on {P}itman-{Y}or processes.} In
ACL-44: \emph{Proceedings of the 21st International Conference on Computational Linguis-
tics and the 44th annual meeting of the Association for Computational Linguistics},
pages 985-992, Morristown, NJ, USA. Association for Computational Linguistics.
\bibitem[{Teh and Jordan(2009)}] {Teh09} Teh Y.W. and Jordan M.I., (2009) {Hierarchical Bayesian nonparametric models with applications.} In N. Hjort, C. Holmes, P. Mueller and S. Walker (Eds.), \emph{Bayesian Nonparametrics: Principles and Practice}, Cambridge, UK: Cambridge University Press, to appear.
\bibitem[{Wallach et al.(2008)}] {Wallach08} Wallach H., Sutton, C. and McCallum, A. (2008) {Bayesian {M}odeling of {D}ependency {T}rees {U}sing {H}ierarchical {P}itman-{Y}or Priors.} In \emph{Proceedings of the Workshop on Prior Knowledge for Text and language (held in conjunction with ICML/UAI/COLT)}, pp. 15--20. Helsinki, Finland, 2008.
\bibitem[{Yau et~al.(2011)}]{Yau11} Yau C., Papaspiliopoulos 0., Roberts G. O., and  Holmes, C. (2011) Bayesian Nonparametric Hidden Markov Models with application to the analysis of copy-number-variation in mammalian genomes, \emph{J. R. Stat. Soc. Series B}, 73(1): 37--57.
\bibitem[{Yau and Holmes(2013)}]{YauHolmes13} Yau C., Holes C.C. (2013) {A decision-theoretic approach for segmental classification}, \textit{Annals of Applied Statistics}, {\bf 7}, 3, 1814--1835.
\bibitem[{Yu(2010)}]{Yu10} Yu S-Z. (2010), Hidden semi-Markov models. \emph{Artificial Intelligence}, 174(2): 215--243.
\end{thebibliography}


\appendix
\section{Appendix: Details of posterior MCMC sampling for the Beta-GOS model} %
\label{app:mcmc}

Here, we provide the details of the MCMC sampling algorithm described in Section \ref{sec:MCMCsampling} for the special case of a Normal sampling distribution and a Normal (or Normal-Inverse-Gamma) base measure.

\subsection{Full conditionals for the Gibbs sampler}
\label{app:gibbscalc}

At each iteration of Gibbs sampler we sample from the full conditionals of $C_n$ and $W_n$, for $n=1,\dots,N$. Here we derive the analytical form of these distributions, for the Beta-GOS model specified in Section \ref{sec:sim}.
Recall that the full conditional distribution for $C_n$ is
\[
\begin{split}
 P\{C_n=i& |C_{-n},W(N),Y(N),  \tau^2\} \\
 &  \propto P\{Y(N) | C_n=i, C_{-n},W(N), \tau^2\} \cdot P\{C_n=i | C_{-n},W(N)\}, \\
\end{split}
\]
where the factor on the right is given by \eqref{eq:C_n} and \eqref{betagos1},
and the left factor is obtained by integration,
\[
\begin{split}
& \hspace{-40pt}P\{Y(N) \mid C_n=i, C_{-n},W(N), \tau^2\} = P\{Y(N) \mid C_n=i, C_{-n}, \tau^2\}\\
&  = \int P\{Y(N),\mu \mid C_n=i, C_{-n}, \tau^2\} ~d\mu \\
& = \prod_{j=1}^{J}\int \prod_{l \ \in \Pi_j} p(Y_l|\mu_j^*) ~ G_0(d\mu_j^*)
\\
& \propto
  \prod_{j=1}^{J} \frac{1}{\sqrt{2\pi}}\frac{1}{( \tau)^{|\Pi_j|}}
  \exp\Big\{-\frac{\sum_{l\in  \Pi_j}y_l^2}{2\tau^2}-\frac{\mu_0^2 }{2\sigma_0^2}+\frac{1}{2}\frac{(\frac{\mu_0}{\sigma_0^2}+
  \sum_{l\in\Pi_j}\frac{{y}_l}{\tau^2})^2}{\frac{1}{\sigma_0^2}+\frac{|\Pi_j|}{\tau^2}}\Big\}
  \frac{1}{\sqrt{{\frac{|\Pi_j|\sigma_0^2}{\tau^2}+{1}}}} ~,
\\
\end{split}
\]
where $\Pi_j$ is the set of indices of data points in cluster $j$, and $J$ is the number of clusters at that iteration.
Note that the latent reinforcements $W(N)$ are used to define the cluster assignments through the data-pairing labels $C(N)$. Conditionally on the data-pairing labels $C(N)$, the data $Y(N)$ is independent of the latent reinforcements $W(N)$.

The full conditional for $W_n$, denoted by $P(W_n|C(N),W_{-n},Y(N))$, is  Beta distributed with updated parameters $A_n,B_n$, defined as in \eqref{eq:C_n}.

\subsection{Inference on the cluster centroids of the Beta-GOS process.}

For the purpose of computational efficiency, it is generally preferable to sample the random partitions integrating out with respect to the parameters of the Beta-GOS process, as described in Section \ref{sec:MCMCsampling} and in Appendix \ref{app:gibbscalc}. If the sampling distribution and the base measure are conjugate, this usually results in improved mixing of the chain. However, in many cases, it may be required to draw inferences on the cluster centroids themselves. As usual with mixtures of DP, inference on the cluster centroids can be easily conducted (even ex-post) from the clustering configurations at each iteration. Therefore, we do not have to sample the centroids within each Gibbs iteration, but if the need be, we can easily resample them at the end of each iteration, or at the end of the sampler from the stored output.

\subsection{Inference on the cluster and global variances}
\label{app:hyper}

Let the variance of the sampling distribution be $\tau^2$.
We assume $\tau^2\sim IGamma(a_0,b_0)$.
The posterior distribution of the variance in each cluster $j$, is given by
\[
 \tau_j^2\mid\mu_j^*, Y_{i} ,i\in \Pi_j \sim IGamma\bigm(a_0 +\frac{|\Pi_j|}{2}, b_0 + \frac{1}{2}\sum_{i \in \Pi_j}(Y_{i}-\mu^*_j)^2 ),
\]

Note that, in case of need and for computational efficiency, we could use these  also quantities to obtain a global estimate for the sampling variance at each iteration, in an MCMC-EM step, as
\(
 \hat{\tau}^2 = \sum_{j=1}^{J}\frac{(|\Pi_j|-1) \tau_j^2}{N- J}.
\)
This may turn useful, for example, for parallelization purposes, as in the simulations of Section \ref{sec:sim}.


\subsection{Inference on the cluster means}
\label{app:atoms}

In the normal-normal model described in Section 5, the posterior distribution of $\mu_j^*$ given data $Y_{i}$ in the $j$-th cluster can be evaluated at each iteration as
\[
 P(\mu_j^* \mid \tau_j^2,Y_{i}, i\in\Pi_j) \sim
N\left(\frac{\frac{\mu_0}{\sigma_0^2}+\frac{|\Pi_j|\overline{Y_j}}{\tau^2_j}}{\frac{1}{\sigma_0^2}+\frac{|\Pi_j|}{\tau^2_j}},\left (    \frac{1}{\sigma_0^2}+\frac{|\Pi_j|}{\tau^2_j}     \right )^{-1}\right) ~.
\]
  for $j=1,\dots,J$, where $\bar{Y}_j$ is the $j$-th cluster specific mean. Note that we have assumed a common sampling variance $\tau^2$; the modification of the previous formula to take into account a  cluster specific variance is of course straightforward.

\section{Appendix: Details of the Proofs and additional theoretical results}

\subsection{Generalized Ottawa Sequence and its moments}

According to \cite{bas-cri-lei}  a sequence $(X_n)_{n \geq 1}$ of random variables taking values in a Polish space
is a Generalized Ottawa Sequence  if there exists a sequence $( W_n)_{n \geq 1}$
(of random variables)
such that the following conditions are satisfied:
{\it {\rm (i)} the law of $X_1$ is $G_0$;  {\rm (ii)} for  $n\geq 1$,
 $X_{n+1}$ and the subsequence $( W_{n+j})_{j \geq 1}$ are conditionally independent given the filtration
 $\CF_n:=\s(X_1,\dots,X_n,W_1,\dots,W_n)$;
{\rm (iii)} the predictive distribution of $X_{n+1}$ given $\CF_n$ is given by
\eqref{betagos1} where the $r_n$'s are  strictly positive functions, $r_n(W_1,\dots,W_n)$, of the vector of latent variables,  such that
 \begin{equation}\label{gos-functions}
 r_n(W_1,\dots,W_n) \geq r_{n+1}(W_1,\dots,W_n,W_{n+1}),
 \end{equation}
 almost surely, with $r_0=1$, and the weights $p_{n,i}=p_{n,i}(W_1,\dots,W_n)$  are
  \begin{equation}\label{gos-functions2}
 \begin{split}
p_{n,i}=\frac{r_n(r_{i-1}-r_{i})}{r_{i-1}r_i}\qquad i=1,\dots,n.\\
\end{split}
 \end{equation}
}

The predictive distribution \eqref{betagos1}-\eqref{formapesi}
corresponds to choice
\( r_{n}(W_1,\dots,W_n)=\prod_{i=1}^n W_i \)
 where $(W_n)_{n \geq 1}$ is a sequence of independent random variables.

We conclude this Section by providing a general result
for the $k$-th moment and for the moment generating function
of the length $K_n$ of a GOS.
Suppose that the sequence $(X_n)_{n \geq 1}$ is a GOS, with $G_0$ diffuse,
and let $U_j=K_j-K_{j-1}$ with $K_0=0$. Then,  $K_n=\sum_{j=1}^n U_j$ and the joint distribution of $U_1, \ldots, U_n$ conditionally on $r_1, \ldots, r_{n-1}$, is
 \[
 P\{U_1=1,\dots,U_n=e_n|r_1,\dots,r_{n-1}\}=\prod_{i=2}^n r_{i-1}^{e_i}(1-r_{i-1})^{1-e_i},
 \]
 for every vector $(e_2,\dots,e_{n})$ in $\{0,1\}^{n-1}$, since $P(U_1=1)=1$ by definition.
Since $K_1=U_1=1$, it follows that, for every $k\geq 1$,
 \[
 P\{K_{n+1}=k+1\}=P \Big\{\sum_{j=2}^{n+1}U_j=k \Big\}=\sum_{\underline e} E \Big [\prod_{i=1}^{n} r_{i-1}^{e_i}
 (1-r_{i-1})^{1-e_i}  \Big ]
 \]
 where the summation is extended over all sequences ${\underline e}=(e_1,\dots,e_n)$ in $\{0,1\}^n$
 such that $\sum_{i=1}^n e_i=k$.
 Moreover, for every $k \geq 1$ and $n \geq 2$, it is easy to see that
 \begin{equation}\label{momentok-1}
 E[(K_{n+1}-1)^k]=E \Big[ \Big(\sum_{j=2}^{n+1}U_j\Big)^k \Big]
 =\sum_{m=1}^{k \wedge n} m! S(k,m) \phi_{n,m}
 \end{equation}
 where $k \wedge n=\min(k,n)$,
  \begin{equation}\label{defphi}
\phi_{n,m} :=\sum_{1 \leq l_1 < l_2 < \dots < l_m \leq n} E[ r_{l_1}\dots r_{l_m}].
  \end{equation}
 and $ S(k,m):=\frac{k!}{m!}\sum_{\{n_i>0: \sum_{i=1}^m n_i=k\}} \frac{1}{n_1! \dots n_{m}!}$ is
the Stirling number of second kind.
Hence, $E[(K_{n+1}-1)^k]$  depends recursively on functions $\phi_{n-1,m}$, $m=1, \ldots, k$.
It may be interesting to note that, using the well known relation between factorial moments
and ordinary moments (see, e.g.,
Example 2.3 in \citealp{Charalambides}), from \eqref{momentok-1}
one gets, for any $m \leq n$,
\begin{equation}\label{fallingmoment}
E[(K_{n+1}-1)_{(m)}]= m! \phi_{n,m}
\end{equation}
where $(t)_{(m)}=t(t-1)\dots(t-m+1)$ is the falling factorial.
Moreover, since
\[
\sum_{k \geq m} (-t)^{k} \frac{S(k,m)}{k!} =\frac{(e^{-t}-1)^m}{m!},
\]
see e.g. Thm. 2.3 in \cite{Charalambides},
 it follows that the moment generating function of $K_{n+1}$ is given by
\begin{equation}\label{momentgenK}
\begin{split}
M_{n+1}(t)&:=E[e^{-t K_{n+1}}]=e^{-t}E[e^{-t (K_{n+1}-1)}] \\
&=e^{-t} \sum_{k \geq 0} \frac{(-t)^{k}}{k!}E[(K_{n+1}-1)^k] = e^{-t} +
e^{-t}\sum_{k \geq 1} \sum_{m=1}^{k \wedge n} \frac{(-t)^{k}m!}{k!} S(k,m) \phi_{n,m} \\
&=e^{-t}+  e^{-t}\sum_{m=1}^{n} m! \,\, \phi_{n,m}\sum_{k \geq m }  \frac{(-t)^{k}}{k!} S(k,m) \\
&=e^{-t}\sum_{m=0}^{n} (e^{-t}-1)^m \phi_{n,m}
\end{split}
\end{equation}
with $\phi_{n,0}:=1$.

\subsection{Proof of Proposition \ref{prop-Kn}}
If we consider equation \eqref{defphi}
with $(W_i)_{i \geq 1}$ independent random variables taking values in $[0,1]$, then
\begin{equation}\label{momentiIND}
\phi_{n,m}= \sum_{1 \leq l_1 < l_2 < \dots < l_m \leq n}
\prod_{j=1}^m \prod _{i=l_{j-1}+1}^{l_j} E[W_{i}^{m+1-j}],
\end{equation}
where $l_0:=0$.
We need some preliminary results.
\begin{lemma}\label{LemmaMOMENTI} If $W_i\sim {\Be}(i+\th-1,1)$, for  given $\th>0$,
then
\begin{equation}\label{formula1}
\phi_{n,m}=\frac{\Gamma(\th+m)}{\Gamma(\th)}
\sum_{j_1=m}^n \sum_{j_2=m}^{j_1} \sum_{j_3=m}^{j_2}
\cdots \sum_{j_m=m}^{j_{m-1}} \frac{1}{(j_1+\th)(j_2+\th) \cdots (j_m+\th)}.
\end{equation}
In particular, as $n$ goes to $+\infty$,
\begin{equation}\label{asmoment}
E[K_n^k]=\frac{\Gamma(\th+k)}{\Gamma(\th)} \log^{k}(n)[1+o(1)].
\end{equation}
\end{lemma}
\noindent

Let us start by proving \eqref{formula1}. First, note that since $W_i$ is a $\Be(i+\th-1,1)$ random variable then, for $1\leq j \leq m$,
$
E[W_{i}^{m+1-j}]=\frac{i+\th-1}{i+\th+m-j}.
$
Hence, by \eqref{momentiIND},
\begin{equation}\label{momentiIND-2}
\phi_{n,m}= \sum_{1 \leq l_1 < l_2 < \dots < l_m \leq n}
\prod_{j=1}^m \prod _{i=l_{j-1}+1}^{l_j} \frac{i+\th-1}{i+\th+m-j}
\end{equation}
which, after some algebra, returns \eqref{formula1}. In order to prove the second part of Lemma \ref{LemmaMOMENTI} we need to introduce additional notation.
For $\th>0$, $k \geq 1$, $m \geq 2$ and $n \geq k$,  set
\[
\begin{split}
\Psi_{k,\th}(n,m)&:=\sum_{j_1=k}^n \sum_{j_2=k}^{j_1} \sum_{j_3=k}^{j_2}
\cdots \sum_{j_m=k}^{j_{m-1}} \frac{m!}{(j_1+\th)(j_2+\th) \cdots (j_m+\th)}, \\
\Psi_{k,\th}(n,1)&:=\sum_{j_1=k}^n \frac{1}{(j_1+\th)}. \\
\end{split}
\]
Note that $m!\phi_{n,m}=\Psi_{m,\theta}(n,m)\Gamma(\theta+m)/\Gamma(\th)$.
For all $k \geq 1$, $m \geq 1$ and $n \geq k$, set
$
Q_{k,\th}(m,n):=\Psi_{k,\th}(n,m)-\log^{m}(n+\th).
$
Now formula \eqref{asmoment} in Lemma \ref{LemmaMOMENTI}
follows easily from \eqref{momentok-1} and the next result.
\begin{lemma}\label{lemmino}
For $\th>0$, $k \geq 1$ and $m \geq 1$, there is a constant $C_{k,\th}(m)$
such that
\begin{equation}\label{asintQ}
|Q_{k,\th}(m,n)| \leq C_{k,\th}(m) \, \log^{m-1}(n+\th) \qquad \text{for every $n\geq k$.}
\end{equation}
\end{lemma}
Let $k \geq 1$ and $\th>0$. For  $m \geq 1$ and $n \geq k$
set
$$
S_{k,\th}(m,n):=\sum_{j=k}^n  \frac{m\log^{m-1}(j+\th)}{j +\th},
$$
and
\begin{equation}\label{re-log1}
R_{k,\th}(m,n):=S_{k,\th}(m,n) -\log^m(n+\th)=\sum_{j=k}^n  \frac{m\log^{m-1}(j+\th)}{j +\th} -\log^m(n+\th).
\end{equation}
We claim that, for any $m \geq 1$, there is a constant $C^*_{m}=C^*_{m,\th,k}$
such that
\begin{equation}\label{re-log2}
|R_{k,\th}(m,n)| \leq C^{*}_m, \qquad \text{for all $n \geq k $}.
\end{equation}
Now observe that
$
\Psi_{k,\th}(n,1)=S_{k,\th}(1,n).
$
Hence, \eqref{re-log2} proves \eqref{asintQ} for $m=1$ and every $k \geq 1$ and $\th>0$.
By induction suppose that \eqref{asintQ} is true for $m=1,\dots,M-1$.
Note that, for $m\geq 2$,
\[
\Psi_{k,\th}(n,m)=\sum_{j_1=k}^n \frac{m}{j_1+\th} \Psi_{k,\th}(j_1,m-1),
\]
hence, by induction hypothesis, for every $\th>0$, $k \geq 1$ and $n \geq k$,
\[
\Psi_{k,\th}(n,M)=\sum_{j_1=k}^n \frac{M}{j_1+\th}
\Big[\log^{M-1}(j_1+\th)+ Q_{k,\th}(M-1,j_1)  \Big].
\]
Using \eqref{re-log1} one gets
\[
\Psi_{k,\th}(n,M)=\log^{M}(n+\th)+R_{k,\th}(M,n)+\sum_{j_1=k}^n \frac{M}{j_1+\th}Q_{k,\th}(M-1,j_1).
\]
Hence, using  \eqref{re-log2} and the induction hypothesis, one can write
\[
\begin{split}
|Q_{k,\th}(M,n)| &\leq |R_{k,\th}(M,n)|+ \sum_{j_1=k}^n \frac{M}{j_1+\th}|Q_{k,\th}(M-1,j_1)|\\
&\leq  C^*_{M,\th,k} + \frac{MC_{k,\th}(M-1)}{M-1} \sum_{j_1=k}^n\frac{M-1}{j_1+\th} \log^{M-2}(j_1+\th)\\
&\leq  C^*_{M,\th,k} + \frac{MC_{k,\th}(M-1)}{M-1}[ \log^{M-1}(n+\th)+ |R_{k,\th}(M-1,n)|  ]\\
&\leq  C^*_{M,\th,k} + \frac{MC_{k,\th}(M-1)}{M-1}[ \log^{M-1}(n+\th)+ C^*_{M-1,\th,k}  ]
\end{split}
\]
which proves \eqref{asintQ} for $m=M$.
To complete the proof let us prove \eqref{re-log2}. Observe that
$
x \mapsto \frac{\log^{m-1}(x+\th)}{x +\th}
$
is a non-increasing function
on $[x_0,+\infty)$ for a suitable $x_0=x_0(k,\th,m)$. Assume, without real loss of generality, that
$k \geq x_0+1$.
Note that, in this case,
\[
\int_{k}^{n+1} \frac{m\log^{m-1}(x+\th)}{x+\th}dx
\leq S_{k,\th}(m,n) \leq \int_{k-1}^{n} \frac{m\log^{m-1}(x+\th)}{x+\th}dx.
\]
Hence,
\[
\log^{m}(n+1+\th)-\log^{m}(k+\th) \leq S_{k,\th}(m,n) \leq \log^{m}(n+\th)-\log^{m}(k-1+\th),
\]
which gives
\[
\log^{m}(n+\th)-\log^{m}(k+\th) \leq S_{k,\th}(m,n) \leq \log^{m}(n+\th),
\]
and then
\[
|S_{k,\th}(m,n) - \log^{m}(n+\th)| \leq \log^{m}(k+\th).
\]

{\it Proof of Proposition \ref{prop-Kn} (a)}. It follows immediately
from  \eqref{asmoment}
and a classical result concerning the convergence in distribution
when the moments converge. 
Indeed,
${E}\left[\left(\frac{K_n}{\log n}\right)^k\right]$ converges to $\frac{\Gamma(\theta+k)}{\Gamma(\theta)}$ that is the $k$-th moment of a $\Gamma(\theta,1)$ random variable.

\vskip 0.5cm

{\it Proof of Proposition \ref{prop-Kn} (b)}.
The first part of the statement of Proposition \ref{prop-Kn}(b) follows from Proposition 2.1
in \cite{bas-cri-lei}
if one shows that
$E[\sum_{i=1}^\infty r_i]<\infty$.
For $\alpha_n=a$ and $\beta_n=b$ one gets
$E[r_n]=a^n/(a+b)^n$ and the thesis follows. When $\alpha_n=n+\theta-1$ and $\beta_n=\beta$,
as explained in Section \ref{s:priorproperties},  $E[r_n]\sim n^{-\beta}$ and the thesis follows since
$\beta>1$.
It remains to prove the assertion concerning the moment generating function
and the factorial moments of $K_\infty$.

 If $\alpha_n=a$ and $\beta_n=b$, \eqref{momentiIND} becomes
\[
\begin{split}
\phi_{n,m}&=\sum_{1 \leq l_1 < l_2 < \dots < l_m \leq n}
\prod_{j=1}^m (E[W_{1}^{m+1-j}])^{l_j-l_{j-1}}, \\
&=\sum_{1 \leq l_1 < l_2 < \dots < l_m \leq n}
\prod_{j=1}^m ( \prod_{i=1}^{m+1-j} \gamma_i )^{l_j-l_{j-1}} \\
\end{split}
\]
since $E[W_1^m]=\prod_{i=1}^{m} \gamma_i$ for $\gamma_i=(a+i-1)/(a+b+i-1)$.
Taking the limit for $n\to +\infty$, we get
\[
\begin{split}
\lim_n \phi_{n,m}& =\sum_{1 \leq l_1 < l_2 < \dots < l_m}
\prod_{j=1}^m \Big ( \prod_{i=1}^{m+1-j} \gamma_i \Big)^{l_j-l_{j-1}}
\\
&= \sum_{k_1\geq 1}\sum_{k_2\geq 1} \cdots  \sum_{k_m \geq 1}
\prod_{j=1}^m  \Big ( \prod_{i=1}^{m+1-j} \gamma_i \Big)^{k_j}\\
&
=\prod_{j=1}^m
\sum_{k_j \geq 1}
 \Big( \prod_{i=1}^{m+1-j} \gamma_i\Big )^{k_j}
=\prod_{j=1}^m \frac{\gamma_1 \cdots \gamma_j}{1-\gamma_1 \cdots \gamma_j} \\
\end{split}
\]
and then
\[
\lim_n \phi_{n,m}=
\prod_{j=1}^m \frac{(a)^{(j)}}{(a+b)^{(j)}-(a)^{(j)}}
\]
where $(t)^{(j)}=t(t+1)\dots(t+j-1)$ is the rising factorial.
Combining this fact with \eqref{momentgenK} it follows that, in this case,
\[
E[e^{-t K_\infty}] =e^{-t} \sum_{m \geq 0} (e^{-t}-1)^m \prod_{j=1}^m \frac{(a)^{(j)}}{(a+b)^{(j)}-(a)^{(j)}}
\]
In addition \eqref{momentok-1}-\eqref{fallingmoment} give
\[
E[\frac{(K_\infty-1)_{m}}{m!}]=  \prod_{j=1}^m \frac{(a)^{(j)}}{(a+b)^{(j)}-(a)^{(j)}}, \qquad
E[(K_\infty-1)^{k}]= \sum_{m=1}^{k} m! S(k,m) \prod_{j=1}^m \frac{(a)^{(j)}}{(a+b)^{(j)}-(a)^{(j)}}
\]

\subsection{Conditionally identity in distribution of the Beta-GOS hierarchical model}
\begin{proposition}\label{propCIDness}
The sequence $(Y_n)_n$ defined by formula \eqref{hiermodel1}-\eqref{hiermodel2} is conditionally identically distributed with respect to the filtration $\mathcal{H}_n=\sigma(Y(n),W(n),\mu(n))$.
\end{proposition}

\begin{proof} Let
$\mathcal{G}_n=\sigma(W(n),\mu(n))$
and $\mathcal{H}_n=\sigma(W(n),\mu(n),Y(n))$.
We have to prove that for every real, bounded and measurable function $g$
\begin{equation}\label{tesicid}
E(g(Y_{n+j})|\mathcal{H}_n)=E(g(Y_{n+1})|\mathcal{H}_n)
\end{equation}
Now, for every $j>0$
\begin{equation}\label{primacosa}
\mathcal{L}(Y_{n+j}\mid
\mathcal{H}_n, \mu_{n+j})= \mathcal{L}(Y_{n+j}\mid \mu_{n+j})=p(\cdot\mid\mu_{n+j})
\end{equation}
and for every $j$ and $n$
\begin{equation}\label{secondacosa}
\mathcal{L}(\mu_{n+j}\mid
\mathcal{H}_n)=\mathcal{L}(\mu_{n+j}\mid
\mathcal{G}_n)
\end{equation}
As already recalled,  $(\mu_n)_n$ is CID with respect to $\mathcal{G}_n=\sigma(W(n),\mu(n))$.
This means that for every  real, bounded and measurable function $f$
\begin{equation}\label{CIDMU}
E(f(\mu_{n+j})|\mathcal{G}_n)=E(f(\mu_{n+1})|\mathcal{G}_n)
\end{equation}
for all $j\geq 1$, see \cite{be-pra-ri}. Thanks to (\ref{secondacosa}), equality (\ref{CIDMU}) also holds with respect the sigma-field $\mathcal{H}_n$. Indeed,
$$E(f(\mu_{n+j})|\mathcal{H}_n)=E(f(\mu_{n+j})|\mathcal{G}_n)=E(f(\mu_{n+1})|\mathcal{G}_n)=E(f(\mu_{n+1})|\mathcal{H}_n)$$
(\ref{primacosa}) implies that
\begin{equation}\label{CIDINT}
E(g(Y_{n+j})|\mathcal{H}_n,\mu_{n+j})=E(g(Y_{n+j})|\mu_{n+j})=\int g(y)p(dy\mid \mu_{n+j})
\end{equation}
(\ref{CIDMU}) and (\ref{CIDINT}) allow to prove the thesis. Indeed,
\begin{equation*}
\begin{split}
E(g(Y_{n+j})|\mathcal{H}_n)&=E(E(g(Y_{n+j})|\mathcal{H}_n,\mu_{n+j}|\mathcal{H}_n)=E(\int g(y)p(dy\mid \mu_{n+j})|\mathcal{H}_n)\\
&=E(\int g(y)p(dy\mid \mu_{n+1})|\mathcal{H}_n)=E(E(g(Y_{n+1})|\mathcal{H}_n,\mu_{n+1}|\mathcal{H}_n)\\
&=E(g(Y_{n+1})|\mathcal{H}_n)\\
\end{split}
\end{equation*}
\end{proof}
\end{document}